\documentclass[twoside,a4paper]{article}

\usepackage{graphicx}
\usepackage{amssymb}
\usepackage{theorem}
\usepackage{amsmath}
\usepackage{subfigure}
\usepackage{enumerate}
\usepackage{color}

\newcommand{\comment}[1]{}

\newtheorem{theorem}{\sc Theorem}[section]
\newtheorem{corollary}{\sc Corollary}[section]
\newtheorem{lemma}{\sc Lemma}[section]
\newtheorem{proposition}{\sc Proposition}[section]
\newtheorem{observation}{\sc Observation}[section]

\newcommand{\R}{{\mathbb R}}

\newcommand{\Q}{{\mathbb Q}}
\newcommand{\Z}{{\mathbb Z}}
\newcommand{\conv}{{\rm conv}}

\newcommand{\intt}{{\rm int}}

\newcommand{\floor}[1]{\left\lfloor #1 \right\rfloor}
\newcommand{\ceil}[1]{\left\lceil #1 \right\rceil}
\newcommand{\area}{{\rm Area}}

\newcommand{\Rf}[1]{R_f(r^1,\dots, r^{#1})}

\title{A Probabilistic Comparison of the Strength of Split,
  Triangle, and Quadrilateral Cuts (extended version)}

\author{Alberto Del Pia$^\ast$ \and Christian Wagner$^\ast$
  \and Robert Weismantel\footnote{Institute for Operations
    Research, ETH Z{\"u}rich, R{\"a}mistrasse 101, 8092
    Z{\"u}rich, Switzerland}}
  
\begin{document}

\maketitle

\begin{abstract}
  We consider mixed integer linear sets defined by two
  equations involving two integer variables and any number
  of non-negative continuous variables.
  The non-trivial valid inequalities of such sets can be
  classified into split, type 1, type 2, type 3, and
  quadrilateral inequalities.
  We use a strength measure of Goemans \cite{goemans} to
  analyze the benefit from adding a non-split inequality on
  top of the split closure.
  Applying a probabilistic model, we show that the
  importance of a type 2 inequality decreases with
  decreasing lattice width, on average.
  Our results suggest that this is also true for type 3 and
  quadrilateral inequalities.
\end{abstract}


\section{Introduction}

We deal with the mixed integer set
$ P_{I} := \{(x,s) \in \Z^{2} \times \R_{+}^{n}: x = f +
  \sum_{j=1}^{n}{r^{j} s_j} \}$
with $f \in \Q^2 \setminus \Z^2$ and
$r^j \in \Q^2 \setminus \{0\}$ for all $j = 1, \dots, n$.
We refer to $f$ as the \emph{root vertex} and to the
vectors $r^j$ as \emph{rays}.
The motivation for analyzing $P_I$ is that it can be
obtained as a relaxation of a general mixed integer linear
set.
Therefore, valid inequalities for $\conv(P_I)$ give rise to
cutting planes for the original mixed integer set.
A thorough investigation of this model can be found in
Andersen et al.~\cite{anlowewo,Strength}.

We consider the set $\Rf{n} := \conv(\{s \in \R_{+}^{n}:
f + \sum_{j=1}^{n}{r^{j} s_j} \in \Z^2 \})$ which is the
projection of $\conv(P_I)$ onto the space of the
$s$-variables.
A closed convex set $B \subseteq \R^2$ with non-empty
interior is \emph{lattice-free} if the interior of $B$ is
disjoint with $\Z^2$.
Furthermore, $B$ is said to be \emph{maximal lattice-free}
if it is not properly contained in another lattice-free
closed convex set.
Any lattice-free polyhedron $B \subseteq
\R^2$ with $f$ in its interior gives rise to a function
$\psi^B : \R^2 \mapsto \R_+$ which is the
Minkowski functional of $B - f$.
We recall that for a convex set $K$ containing the origin in
its interior, the \emph{Minkowski functional}
$\| \cdot \|_K$ of $K$ is defined by $\| r \|_K := \min
\{\lambda > 0 : r \in \lambda K \}$ where $r \in \R^2$.
We call $\sum_{j=1}^n{\psi^B(r^j)s_j} \ge 1$ the
\emph{cut associated with $B$}.
This inequality is valid for $\Rf{n}$.
Conversely, any non-trivial valid inequality for $\Rf{n}$
can be obtained from a lattice-free polyhedron.

Every lattice-free rational polyhedron $B
\subseteq \R^2$ with $f \in \intt(B)$ is contained in a
maximal lattice-free rational polyhedron $\bar{B}$.
It follows that $\psi^{\bar{B}}(r^j) \le \psi^B(r^j)$ for
all $j = 1, \dots, n$ and the cut associated with $B$ is
dominated by the cut associated with $\bar{B}$.
This implies that the benefit from adding the latter cut on
top of the split closure is at least as high as the benefit
from adding the former.
In this paper, we focus exclusively on cuts associated with
maximal lattice-free rational polyhedra.

A classification of planar maximal lattice-free closed
convex sets was given by Lov\'asz \cite{lovasz}.
The refined classification which is stated below can be
found in \cite{DeyWolsey}.

\begin{figure}[ht]
  \unitlength=1mm
  \begin{center}
    \subfigure[\label{mlf.split} Split]{
        \begin{picture}(23,20)
          \put(2,0){\includegraphics[width=20\unitlength]{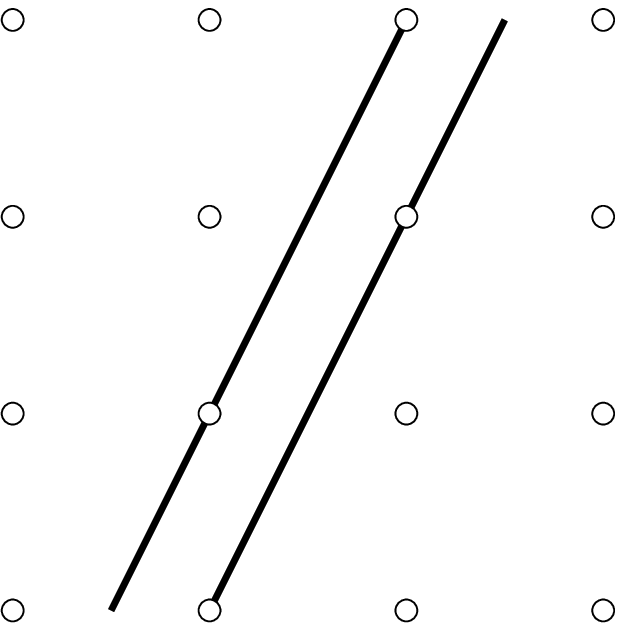}}
        \end{picture}
        } 
        \quad
    \subfigure[\label{mlf.tr1} Type 1 triangle]{
        \begin{picture}(23,20)
          \put(2,0){\includegraphics[width=20\unitlength]{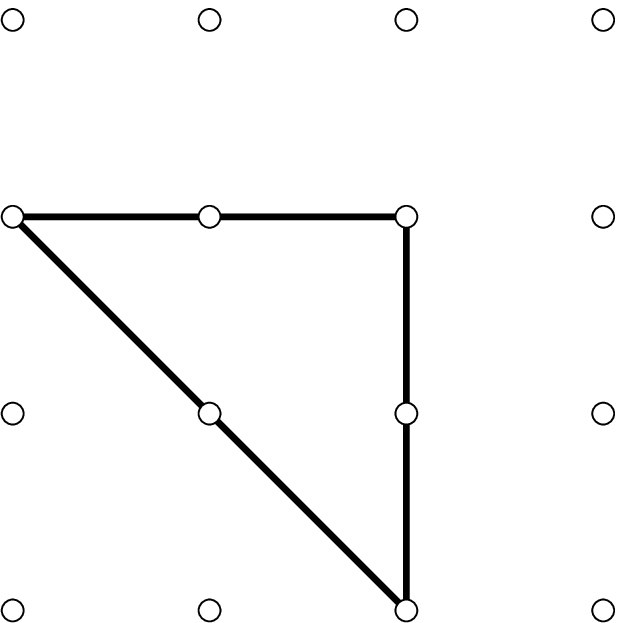}}
        \end{picture}
        } 
    \quad 
    \subfigure[\label{mlf.tr2} Type 2 triangle]{
        \begin{picture}(23,20)
          \put(2,0){\includegraphics[width=20\unitlength]{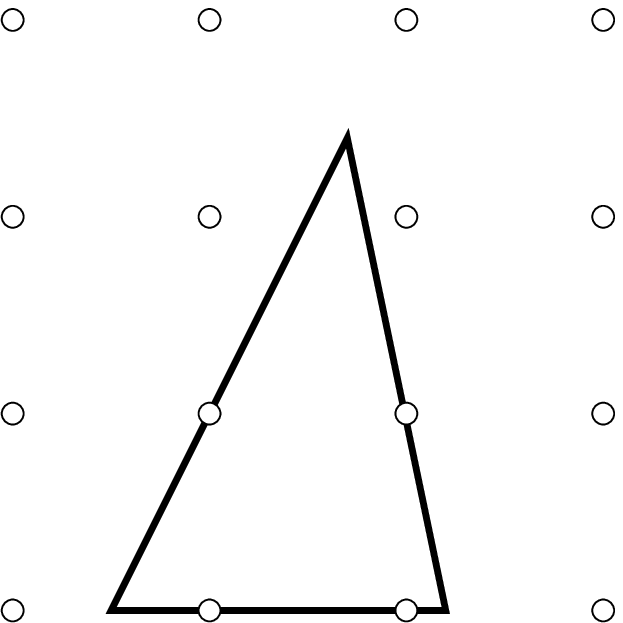}}
        \end{picture}
        }
    \quad
    \subfigure[\label{mlf.tr3} Type 3 triangle]{
        \begin{picture}(23,20)
          \put(2,0){\includegraphics[width=20\unitlength]{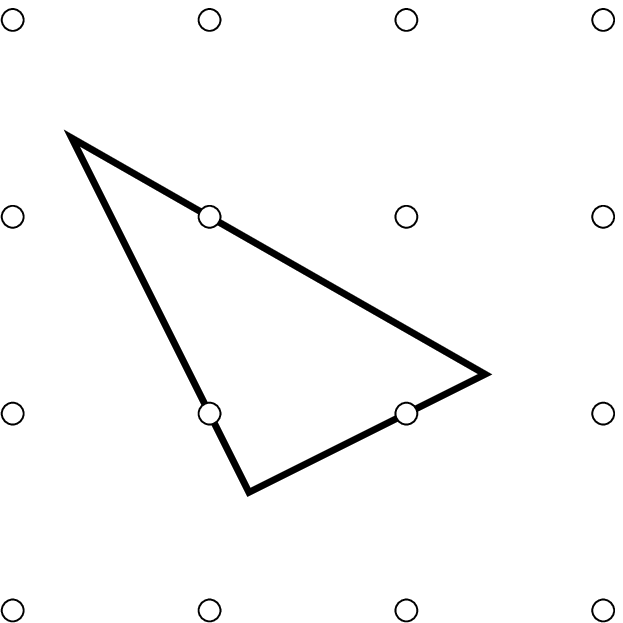}}
        \end{picture}
        } 
    \quad 
    \subfigure[\label{mlf.quad} Quadrilateral]{
        \begin{picture}(23,20)
          \put(2,0){\includegraphics[width=20\unitlength]{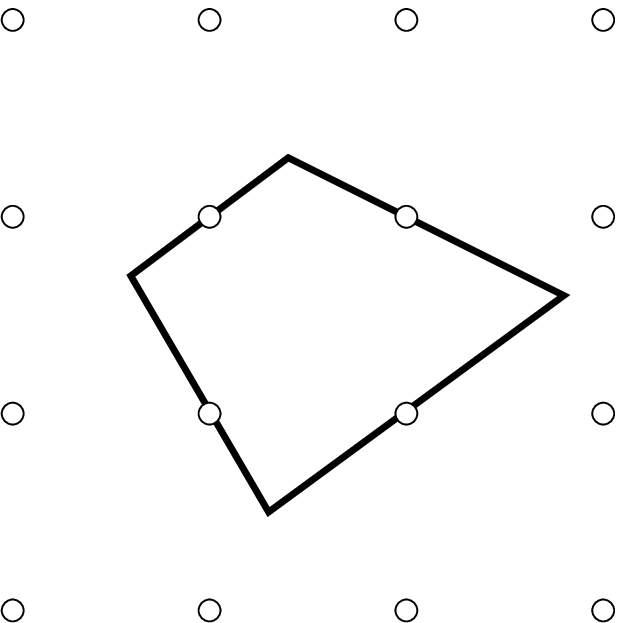}}
        \end{picture}
        }
  \end{center}
  \caption{\label{mlf.sets} All types of maximal
    lattice-free convex sets in $\R^2$ with non-empty interior.}
\end{figure}

\begin{proposition} [\cite{DeyWolsey,lovasz}]
  \label{list-of-two-dim-max}
Let $B$ be a maximal lattice-free closed convex set with
non-empty interior.
Then $B$ is a polyhedron and the relative interior of each
facet of $B$ contains at least one integer point.
In particular, $B$ is one of the following sets (see
Fig.~\ref{mlf.sets}).
\begin{enumerate} [I.]
  \item A split $\{(x_1,x_2) \in \R^2 : b \leq a_1 x_1 + a_2
    x_2 \le b + 1\}$ where $a_1$ and $a_2$ are coprime
    integers and $b$ is an integer.
  \item A triangle which in turn is either
    \begin{enumerate} [(a)]
      \item a type 1 triangle, i.e.~a triangle with integer
        vertices and exactly one integer point in the
        relative interior of each edge, or
      \item a type 2 triangle, i.e.~a triangle with at
        least one fractional vertex $v$, exactly one integer
        point in the relative interior of the two edges
        incident to $v$ and at least two integer points on
        the third edge, or 
      \item a type 3 triangle, i.e.~a triangle with exactly
        three integer points on the boundary, one in the
        relative interior of each edge.
    \end{enumerate}
  \item A quadrilateral containing exactly one integer point
    in the relative interior of each of its edges.
\end{enumerate}
\end{proposition}

By Proposition \ref{list-of-two-dim-max}, it follows that
$\Rf{n}$ has three types of non-trivial valid inequalities:
split, triangle, and quadrilateral inequalities named after
the corresponding two-dimensional object from which the
inequality can be derived as an intersection cut
\cite{balasint}.
Triangle inequalities are further subdivided into type 1,
type 2, and type 3 inequalities (see Fig.~\ref{mlf.sets}).
Note that a non-trivial valid inequality for $\Rf{n}$ can
correspond to more than one maximal lattice-free
polyhedron.

We are interested in the quality of cuts associated with
the maximal lattice-free polyhedra described in
Proposition \ref{list-of-two-dim-max}.
We use the strength measure introduced by Goemans
\cite{goemans} to evaluate the quality of a cut.
Basu et al.~\cite{BasuBonamiCornuejolsMargot09} assess
the strength of split, triangle, and quadrilateral
inequalities in a non-probabilistic setting.
They show that split and type 1 inequalities may produce an
arbitrarily bad approximation of $\Rf{n}$.
On the other hand, type 2 or type 3 or quadrilateral
inequalities deliver good approximations of $\Rf{n}$ in
terms of the strength.
This, however comes with a price.
Up to transformations preserving the integer lattice, there
is only one split and only one triangle of type 1, but an
infinite number of triangles of types 2 and 3 and
quadrilaterals.
Therefore, it can be expected that for real instances an
approximation by adding all these cuts is hard to compute.
From a more practical point of view, one is interested in
approximations of the mixed integer hull that one can
generate easily.
Current state-of-the-art in computational integer
programming is to experiment with split cuts and the split
closure (see e.g.~\cite{AndersenWeismantel,BalasSaxena}).
This is the point of departure of our theoretical study.

The aim of this paper is to shed some light on the question
which average improvement a non-split cut gives when added
on top of the split closure.
For that, we take any maximal lattice-free triangle or
quadrilateral $B$, and we investigate all potential sets
$\Rf{n}$ such that we can generate a valid cut from $B$.
For this, it is required that $f$ is in the interior of
$B$.
We vary $f$ uniformly at random in the interior of $B$.
This is our probability distribution.
For each particular $B$ and $f \in \intt(B)$ we let $n$ and
$r^1, \dots, r^n$ attain arbitrary values.
We compute a lower bound on the probability that the
strength of adding the cut associated with $B$ on top of the
split closure is less than or equal to an arbitrary value. 

As a conclusion from our probabilistic analysis we obtain
that the addition of a single type 2, type 3, or
quadrilateral inequality to the split closure becomes less
likely to be beneficial the closer the lattice-free set
looks like a split, i.e.~the closer its geometric width is
to one.

We think that this complements nicely the analysis in
\cite{BasuBonamiCornuejolsMargot09}:
they construct sequences of examples in
which cuts from triangles of types 2 and 3, and
quadrilaterals cannot be approximated to within a constant
factor by the split closure.
The approximation becomes worse as the triangles and
quadrilaterals converge towards a split.
From the results in our paper it follows that this
geometrically counterintuitive situation occurs extremely
rarely.


\section{Probabilistic model and main results}

Often, the quality of cuts is measured according to their
worst case performance.
In this paper we apply a stochastic approach by introducing
a probability distribution on all possible instances.
See \cite{BasuCornuejolsMolinaro10} and
\cite{HeAhmedNemhauser} for other approaches on how to
apply a probabilistic analysis to mixed integer linear
sets.

Our aim is to evaluate the benefit from adding a single
cut associated with a maximal lattice-free rational triangle
or quadrilateral on top of the split closure.
For that reason we apply the following strength measure.

A non-empty polyhedron $P = \left\{ x \in \R_{+}^{n}: Ax
\geq b \right\}$ is said to be \emph{of covering type} if
all entries in $A$ and $b$ are non-negative and $0 \not \in
P$.
Let $\alpha \in (0,+\infty)$.
We call $\alpha P := \{x: \alpha x \in P \}$ the
\emph{dilation of $P$ by the factor $\alpha$} and define
$\alpha P := \R^n_+$ if $\alpha =+\infty$.
Note that $P \subseteq \alpha P$ for $\alpha \geq 1$.
Let $Q \subseteq \R^n_+$ be any convex set such that $P
\subseteq Q$.
The \emph{strength of $P$ with respect to $Q$}, denoted by
$t(P,Q)$, is the minimum value of $\alpha \geq 1$ such that
$Q \subseteq \alpha P$.

Here, the underlying idea is that a high strength means high
benefit.
In other words, let $P$ and $Q$ be two relaxations for
a mixed integer set $M$ such that $P \subseteq Q$.
The larger $t(P,Q)$ the better is the gap to $M$ closed by
$P$ compared to that of $Q$.

\begin{lemma}[\cite{goemans}] \label{CovType}
Let $P := \{x \in \R_{+}^{n}: a_{i}^{T}x \geq b_{i}\ {\rm for\
all\ } i = 1, \dots, m\}$ be a polyhedron of covering type
and $Q \subseteq \R^n_+$ any convex set such that $P
\subseteq Q$.
Then,
$$ t(P,Q) = \max_{i=1,\dots,m}\left\{\frac{b_{i}}{\inf\{a_i^T x : x
   \in Q\}} : b_i > 0\right\}. $$
If $\inf\{a_i^T x : x \in Q\} = 0$ for some $i \in
\{1,\dots,m\}$, where $b_i > 0$, then $t(P,Q)$ is defined to
be $ + \infty$.
\end{lemma}

In the following we use $R_f^n$ instead of $\Rf{n}$ for
simplicity.
Note that $R_f^n$ is a polyhedron \cite{Meyer74}.
Let $\Omega$ be the set of all maximal lattice-free rational
polyhedra in $\R^2$ which contain $f$ in the interior and
let ${\cal S}$ (resp.~${\cal T}_1$, ${\cal T}_2$,
${\cal T}_3$, ${\cal Q}$) be the subset of $\Omega$
containing all splits (resp.~triangles of type 1, type 2,
type 3, quadrilaterals).
Observe that $R_f^n = \{ s \in \R^n_+ :
\sum_{j=1}^n{\psi^B(r^j) s_j} \ge 1 \textup{ for all } B \in
{\cal S} \cup {\cal T}_1 \cup {\cal T}_2 \cup {\cal T}_3
\cup {\cal Q} \}$.
For a non-empty set ${\cal L} \subseteq \Omega$ we
define ${\cal L}(R_f^n)$ to be the intersection of all
cuts associated with the polyhedra in ${\cal L}$ and the
trivial inequalities $s_j \ge 0$ for all $j = 1, \dots, n$.
In the remainder of the paper, ${\cal L}$ will always be
${\cal S}$ or ${\cal S} \cup \{B\}$ for some
$B \in \Omega \setminus {\cal S}$.
In this case ${\cal L}(R_f^n)$ is a polyhedron
\cite{CooKanSch90} of covering type.

In order to evaluate the gain from adding a triangle or
quadrilateral inequality to the split closure we apply the
following procedure.
We fix a maximal lattice-free rational triangle or
quadrilateral $B \in \Omega \setminus {\cal S}$, but we
allow $f$ to vary in its interior.
We will show that by considering a specific set of rays
which depends only on $f$, we obtain an upper bound on
$t({\cal F}(R_f^n), {\cal S}(R_f^n))$ where ${\cal F} =
{\cal S} \cup \{B\}$.
Depending on where $f$ is located this bound may differ.
By varying $f$ over the entire area of $B$ we compute the
area for which the bound is below a certain value, $z$
say, and compare it to the area of $B$.
This gives a ratio which is a lower bound on the probability
that the strength is less than or equal to $z$.
In turn, $1$ minus this probability is an upper bound for
the chance that $B$ improves upon the split closure by a
value of more than $z$ with respect to the strength.

Let $B \in \Omega \setminus {\cal S}$ and let ${\cal F} =
{\cal S} \cup \{B\}$.
The following observation follows easily from Lemma
\ref{CovType}.

\begin{observation} \label{obs}
$$t({\cal F}(R_f^n), {\cal S}(R_f^n))
= \frac{1}{\min \{ \sum_{j=1}^n{\psi^B(r^j) s_j}:
s \in {\cal S}(R_f^n) \} }.$$
\end{observation}

Since $B$, $f$, and $r^1, \dots, r^n$ are rational, we can
assume that the rays $r^1, \dots, r^n$ are scaled such that
the points $f + r^j$ are on the boundary of $B$ (see
\cite{BasuBonamiCornuejolsMargot09} for an explanation why
this assumption is feasible).
In the course of the paper we will deal with optimization
problems of the following type:
\begin{equation} \label{inf-prob}
  \min\ \sum_{j=1}^n{\psi^B(r^j) s_j
  \textup{\quad s.t.\quad} s \in {\cal S}(R_f^n)}
\end{equation}
for some $B \in \Omega \setminus {\cal S}$.
Scaling the rays as described above implies that the
objective function becomes $\sum_{j=1}^n{s_j}$.

We define a \emph{corner ray} to be a ray $r^j$ where the
point $f + r^j$ is a vertex of $B$.
For now, assume that all the corner rays (and only those)
are present.
We will explain later why we can make this assumption.
This implies that the corner rays are fixed once $f$ and $B$
are chosen.
We assume a continuous uniform distribution on $f$ in the
interior of $B$.
Given $z \in \R$, $z > 1$, we define $P^B(z)$ to be the
probability that 
$t({\cal F}(R_f^n), {\cal S}(R_f^n))$ is less than or
equal to $z$ for $f$ varying in the triangle or
quadrilateral $B$, i.e.
$$ P^B(z) := \frac{1}{\area(B)}
   \int\limits_{f \in \intt(B)}{\boldsymbol{1}\{t({\cal
   F}(R_f^n), {\cal S}(R_f^n)) \le z\}df},$$
where $\area(B)$ is the area of $B$ and $\boldsymbol{1}$ is
the indicator function.

Let us now argue why we can assume that all corner rays
of $B$ are present in order to obtain the desired bound.
Assume $k \ge 1$ corner rays $r^{n+1}, \dots, r^{n+k}$ are
missing and consider $R_f^{n+k} := \conv(\{s \in
\R_{+}^{n+k}: f + \sum_{j=1}^{n+k}{r^{j} s_j} \in \Z^2 \})$.
We now apply Observation \ref{obs} to infer
$t({\cal F}(R_f^{n+k}), {\cal S}(R_f^{n+k})) \ge 
t({\cal F}(R_f^{n}), {\cal S}(R_f^{n})) \Longleftrightarrow
\min \{ \sum_{j=1}^{n+k}{s_j}: s \in
{\cal S}(R_f^{n+k}) \} \le
\min \{ \sum_{j=1}^n{s_j}: s \in {\cal S}(R_f^n) \}$.
The latter inequality follows from the fact that an optimal
solution $\bar{s}$ for the latter minimization problem
implies a feasible solution for the former minimization
problem by setting $\bar{s}_j = 0$ for all $j = n+1, \dots,
n+k$.
This is because if there would be a split which cuts off
$(\bar{s},0)$ from ${\cal S}(R_f^{n+k})$, then the same
split would cut off $\bar{s}$ from ${\cal S}(R_f^{n})$.
It follows that
$\boldsymbol{1} \{ t({\cal F}(R_f^{n+k}),
{\cal S}(R_f^{n+k})) \le z \} \le
\boldsymbol{1} \{ t({\cal F}(R_f^{n}),
{\cal S}(R_f^{n})) \le z \}$.
Thus, by adding the corner rays we obtain a lower bound on
$P^B(z)$.

We now explain why we can assume that only the corner rays
of $B$ are present in order to compute
$t({\cal F}(R_f^{n+k}), {\cal S}(R_f^{n+k}))$.
In \cite{BasuBonamiCornuejolsMargot09} it is shown that, if
all corner rays are present and under the assumption of
scaled rays, \eqref{inf-prob} reduces to the problem where
only the corner rays are given, i.e.~the
objective function is $\sum_{j=1}^{k}{s_j}$, where
$\{r^1,\dots,r^k\}$, $k \in \{3,4\}$ is exactly the set of
corner rays for the given triangle or quadrilateral $B$.
By Observation \ref{obs} this implies
$t({\cal F}(R_f^{n+k}), {\cal S}(R_f^{n+k})) =
t({\cal F}(R_f^{k}), {\cal S}(R_f^{k}))$.
Thus, in the remainder of the paper, in order to obtain the
desired bounds, we assume that the set of rays consists of
exactly the corner rays of $B$.

An important value in our analysis is the so-called
\emph{lattice width}.
Let $B \subseteq \R^2$ be a maximal lattice-free set with
non-empty interior.
The lattice width of $B$ is defined to be
\begin{equation} \label{lwidth}
  w(B) := \min_{u \in \Z^2 \setminus \{0\}}\left\{
  \max_{x \in B}u^Tx - \min_{x \in B}u^Tx \right\}.
\end{equation}
It is shown in \cite{Hurkens} that $1 \le w(B) \le 1 +
\frac{2}{\sqrt{3}}$.
Furthermore, applying an affine unimodular transformation
(i.e.~a 
transformation preserving $\Z^2$) to $B$ changes neither the
lattice width nor the strength.
In the remainder of the paper, we informally call $B$
\emph{flat} whenever $w(B)$ is sufficiently close to $1$.
Our main results are summarized below.

\begin{theorem}[probabilistic strength -- lattice width
  relation] \label{main.th}
Let $T_i$ be a triangle of type $i \in \{1,2\}$ and $w :=
w(T_2)$.
Then, for any $z > 1$, we have
\begin{enumerate}[I.]
  \item \label{th.1} \begin{equation*}
          P^{T_1}(z) =
          \begin{cases}
            0 & \textup{if } 1 < z \le \frac{3}{2}, \\
            \frac{3}{4} (\frac{2z-3}{z-1})^2 & \textup{if }
            \frac{3}{2} < z < 2, \\
            1 & \textup{if } 2 \le z < +\infty. \\
          \end{cases}
        \end{equation*}
  \item \label{th.t2} \begin{equation*}
          P^{T_2}(z) \ge
          \begin{cases}
            0 & \textup{if } 1 < z \le w, \\[1.8mm]
            g_1 &
            \textup{if } w < z \le \frac{w}{w-1}, \\[1.8mm]
            g_1 + g_2 &
            \textup{if } \frac{w}{w-1} < z < + \infty, \\
          \end{cases}
        \end{equation*}
with $g_1 = \frac{(z-w)(2wz-w-z)}{w^2(z-1)^2}$ and
$g_2 = \frac{(w-1)^2(z-1)^2-1}{w^2(z-1)^2}$.
\end{enumerate}
\end{theorem}

Theorem \ref{main.th} \ref{th.t2} shows that for any given
$z > 1$, $P^{T_2}(z)$ tends to $1$ if $w(T_2)$ converges to
$1$, i.e.~the probability that a flat type 2 triangle
improves upon the split closure by a value of more than $z$
goes to $0$.
This will be explained in further detail in Section
\ref{sec.t2}.

The analysis of a type 3 triangle $T_3$ and a quadrilateral
$Q$ turns out to be more complex.
We did not succeed in putting $P^{T_3}(z)$
(resp.~$P^{Q}(z)$) into direct relation to $w(T_3)$
(resp.~$w(Q)$) and $z$ only.
Instead we parametrize $T_3$ and $Q$ in terms of the
coordinates of their vertices. 
Using this more complicated parametrization we derive
formulas for $P^{T_3}(z)$ and $P^{Q}(z)$.
Then we discretize the coordinates of the vertices and
evaluate the formulas with respect to our discretization.
This qualitatively leads to the same conclusion as before:
if $T_3$ and $Q$ converge towards a split (meaning the
lattice width converges to $1$), the probability
that the strength of the associated inequality is less than
or equal to $z$ tends to $1$.
We refer to Sections \ref{sec.quad} and \ref{sec.t3} for
the corresponding formulas.


\section{Type 1 triangles}

By an affine unimodular transformation,
we assume that the type 1 triangle $T_1$ is given by
$T_1 = \conv\{(0,0),(2,0),(0,2)\}$.
Let $R_1 := \intt(\conv \{(1,0),(0,1),(1,1)\})$,
$R_2 := \intt(\conv \{(0,0),(1,0),(0,1)\})$,
$R_3 := \intt(\conv \{(0,1),(1,1),(0,2)\})$, and
$R_4 := \intt(\conv \{(1,0),(1,1),(2,0)\})$.

Note that $\intt(T_1) \setminus \cup_{j=1}^{4}{R_j}$ is a
set of area zero, so it can be neglected.
For given $f \in R_1 \cup R_2 \cup R_3 \cup R_4$ the three
corner rays are $r^1 = (-f_1,-f_2)$, $r^2 = (2-f_1,-f_2)$,
and $r^3 = (-f_1,2-f_2)$.
Let ${\cal F} = {\cal S} \cup \{T_1\}$. 
In \cite{BasuBonamiCornuejolsMargot09} it is shown that
\begin{equation} \label{t1.strength}
   t({\cal F}(R_f^3), {\cal S}(R_f^3)) =
  \begin{cases}
    2 & \textup{if } f \in R_1, \\
    \frac{3-f_1-f_2}{2-f_1-f_2} & \textup{if } f \in R_2, \\
    \frac{f_2+1}{f_2} & \textup{if } f \in R_3, \\
    \frac{f_1+1}{f_1} & \textup{if } f \in R_4. \\
  \end{cases}
\end{equation}
For the sake of a shorter notation we write
$\boldsymbol{1}\{t \le z\}$ instead of
$\boldsymbol{1}\{t({\cal F}(R_f^3), {\cal S}(R_f^3)) \le
z\}$.
Then,
$$ P^{T_1}(z) = \frac{1}{\area(T_1)}
   \sum_{j=1}^{4}{\int\limits_{f \in R_j}{\boldsymbol{1}\{t
   \le z\}df}}. $$ 
We compute the four integrals separately.
For that we need to check when the corresponding functions
in \eqref{t1.strength} attain a value which is less than or
equal to $z$.
Assume $f \in R_1$.
Then $\int_{f \in R_1}{\boldsymbol{1}\{t \le z\}df}$ is $0$
if $z < 2$ and is $\frac{1}{2}$ if $z \ge 2$.
Assume $f \in R_2$.
We have $\frac{3-f_1-f_2}{2-f_1-f_2} \le z \Leftrightarrow
f_1 + f_2 \le \frac{2z-3}{z-1}$
and thus $\int_{f \in R_2}{\boldsymbol{1}\{t \le z\}df}$ is
$0$ if $z \le \frac{3}{2}$, is $\frac{1}{2}$ if $z \ge 2$,
and is $\frac{1}{2} (\frac{2z-3}{z-1})^2$ otherwise (here we
used the fact that $0 < f_1 + f_2 < 1$ in $R_2$).
Assume $f \in R_3$.
Then $\frac{f_2+1}{f_2} \le z \Leftrightarrow f_2 \ge
\frac{1}{z-1}$.
Hence, $\int_{f \in R_3}{\boldsymbol{1}\{t \le z\}df}$ is
$0$ if $z \le \frac{3}{2}$, is $\frac{1}{2}$ if $z \ge 2$,
and is $\frac{1}{2} (\frac{2z-3}{z-1})^2$ otherwise (here we
used $1 < f_2 < 2$).
Finally, the case $f \in R_4$ is analogous to the previous
case.
Since $\area(T_1) = 2$ it follows
\begin{equation*}
   P^{T_1}(z) =
  \begin{cases}
    0 & \textup{if } 1 < z \le \frac{3}{2}, \\
    \frac{3}{4} (\frac{2z-3}{z-1})^2 & \textup{if }
    \frac{3}{2} < z < 2, \\
    1 & \textup{if } 2 \le z < +\infty. \\
  \end{cases}
\end{equation*}


\section{Proof strategy}

The analysis of the strength of type 1 triangles is quite
easy since the
split closure is known: it is always defined by a subset of
the three split inequalities corresponding to the splits
whose normal vectors are the normal vectors of the facets of
the type 1 triangle (see \cite{BasuBonamiCornuejolsMargot09}
for a proof).
Note that $w(T_1) = 2$ and that this value is attained by
precisely these vectors.

Using the split closure for triangles of types 2 and 3,
and quadrilaterals would result in too complicated
formulas.
Thus, we choose another strategy.
Instead of using the entire split closure we will take only
one well-chosen split inequality and therefore
obtain lower bounds for the desired probabilities.
Let $B \in \Omega \setminus {\cal S}$.
The split inequality which we choose will depend on the
location of $f$ in the interior of $B$.
For that we partition $B$ into regions $R_1, \dots, R_p$ and
select a single split for each region.
The basic idea is to choose a split that contains $f$ in its
interior and such that the normal vector of the split is a
potential candidate for an integer vector for which $w(B)$
is attained.
Our candidate splits will always be among the vectors
$(1,0)$, $(0,1)$, and $(1,1)$ since we apply an affine
unimodular transformation to $B$ to bring it into an
appropriate form (see Sections \ref{sec.t2}, \ref{sec.quad},
and \ref{sec.t3} for details).

We now show that our simplification of using only one split
inequality instead of the split closure leads to a lower
bound for $P^B(z)$.
Let ${\cal F} = {\cal S} \cup \{B\}$ for some $B \in \Omega
\setminus {\cal S}$ and let $L \in {\cal S}$ be arbitrary.
Then, by Observation \ref{obs},
$t({\cal F}(R_f^k), {\cal S}(R_f^k))
= ( \min \{ \sum_{j=1}^{k}{s_j}:
s \in {\cal S}(R_f^k) \} )^{-1} \le
( \min \{ \sum_{j=1}^{k}{s_j}:
s_j \ge 0 \textup{ for } j \in \{1, \dots, k\}
\textup{ and } \sum_{j=1}^k{\psi^L(r^j)s_j} \ge 1 \})^{-1} 
= t(\{ \sum_{j=1}^k{\psi^L(r^j)s_j} \ge 1,
  \sum_{j=1}^k{\psi^B(r^j)s_j} \ge 1 \} \cap \{s_j \ge 0
  \textup{ for } j \in \{1, \dots, k\}\},
  \{\sum_{j=1}^k{\psi^L(r^j)s_j} \ge 1 \} \cap \{s_j \ge 0
  \textup{ for } j \in \{1, \dots, k\} \})$,
since we assumed $\psi^B(r^j) = 1$ for all $j = 1, \dots, k$.

For ease of notation we denote the latter by $\bar{t}(B,L)$,
which is the strength of the polyhedron obtained by adding
the cuts associated with $B$ and $L$ with respect to the
polyhedron obtained by just adding the cut associated with
$L$.
It follows $\boldsymbol{1}\{t({\cal F}(R_f^k),
{\cal S}(R_f^k)) \le z\} \ge \boldsymbol{1}\{\bar{t}(B,L) \le
z\}$ and therefore 
\begin{equation} \label{T2.approx}
  P^{B}(z) \ge \frac{1}{\area(B)} \sum_{j=1}^{p}{
  \int\limits_{f \in R_j}{\boldsymbol{1}\{\bar{t}(B,L_{R_j}) \le z\}df}},
\end{equation}
where $L_{R_j}$ is the single split which is used in region
$R_j$ to approximate the split closure, for
$j = 1, \dots, p$.
In the following, for simplicity, we write
$\boldsymbol{1}\{\bar{t} \le z\}$ instead of
$\boldsymbol{1}\{\bar{t}(B,L) \le z\}$ whenever $B$ and $L$ are
clear from the context.
In order to compute $\bar{t}(B,L)$ we need to
solve an optimization problem of the type
\begin{equation*}
  \begin{split}
    \min\ & s_1 + \dots + s_k \\ \mbox{\quad s.t. } &
    \psi^{L}(r^1)s_1 + \dots + \psi^{L}(r^k)s_k  \ge 1, \\
    & s_j \ge 0 \quad \mbox{for } j = 1, \dots, k.
  \end{split}
\end{equation*}
The values $\psi^{L}(r^1), \dots, \psi^{L}(r^k)$ are called
the \emph{coefficients of the split $L$}.
In general, for a split $L = \{x \in \R^2 : \floor{\pi^Tf}
\le \pi^T x \le \ceil{\pi^Tf}\}$, $\pi \in \Z^2 \setminus
\{0\}$, containing $f$ in its interior and with
corresponding inequality $\sum_{j=1}^k \psi^{L}(r^j)s_j \ge
1$, it is easy to verify that
\begin{equation} \label{coefficients}
  \psi^{L}(r^j) =
    \begin{cases}
      \tfrac{\pi^T r^j}{\ceil{\pi^Tf} - \pi^Tf} &
      \textup{if } \pi^T r^j > 0, \\
      0 & \textup{if } \pi^T r^j = 0, \\
      \tfrac{\pi^T r^j}{\floor{\pi^Tf} - \pi^Tf} &
      \textup{if } \pi^T r^j < 0.
    \end{cases}
\end{equation}
Thus, we just need the normal vector $\pi$ of the split $L$
to compute the split coefficients $\psi^{L}(r^1), \dots,
\psi^{L}(r^k)$.

In Section \ref{sec.t2} we will explain our computations for
triangles of type 2 in detail.
Since the computations for quadrilaterals and triangles of
type 3 in Sections \ref{sec.quad} and \ref{sec.t3} give no new
insights we will only state intermediate results there.


\section{Type 2 triangles} \label{sec.t2}

\begin{figure}[ht]
        \centerline{\includegraphics[scale = 1.0]{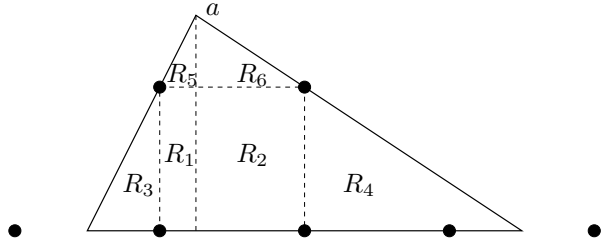}}
        \caption{Decomposition of a type 2 triangle.}
        \label{type2}
\end{figure}

By an affine unimodular transformation, we assume that the
type 2 triangle $T_2$ has one facet containing the points
$(0,0)$ and $(1,0)$, one facet containing $(0,1)$, and one
facet containing $(1,1)$.
Furthermore, one vertex $a = (a_1,a_2)$ satisfies $0 < a_1 <
1$ and $1 < a_2$.
Thus, the other vertices are
$(-\frac{a_1}{a_2-1},0)$ and
$(\frac{a_2-a_1}{a_2-1},0)$.
We assume that $a$ is arbitrary but fixed and treat it as a
parameter in the subsequent computations.
We decompose $T_2$ into six regions:
$R_1 := \intt(\conv\{(0,0),(a_1,0),(a_1,1),(0,1)\})$,
$R_2 := \intt(\conv\{(a_1,0),(1,0),(1,1),(a_1,1)\})$,
$R_3 := \intt(\conv\{(0,0),(0,1),(-\frac{a_1}{a_2-1},0)\})$,
$R_4 := \intt(\conv\{(1,0),$ $(1,1),(\frac{a_2-a_1}{a_2-1},0)\})$,
$R_5 := \intt(\conv\{(0,1),(a_1,1),(a_1,a_2)\})$, and
$R_6 := \intt(\conv\{(a_1,1),(1,1),(a_1,a_2)\})$ (see
Fig.~\ref{type2}).

For given $f \in \cup_{j=1}^{6}{R_j}$ the three corner rays
are $r^1 = (-\frac{a_1}{a_2-1}-f_1,-f_2)$,
$r^2 = (\frac{a_2-a_1}{a_2-1}-f_1,-f_2)$,
and $r^3 = (a_1-f_1,a_2-f_2)$.
Furthermore, $\area(T_2) = \frac{(a_2)^2}{2(a_2-1)}$.
Our assumptions on $T_2$ imply that $w(T_2)$, as defined in
\eqref{lwidth}, is attained for $u = (1,0)$ or $u = (0,1)$
(see \cite{AverkovWagner}).
Thus, $w := w(T_2) = \min\{a_2, \frac{a_2}{a_2-1}\}$ which
implies $w = a_2$ if $a_2 \le 2$ and $w = \frac{a_2}{a_2-1}$
if $a_2 > 2$.
In either case we have $\area(T_2) = \frac{w^2}{2(w-1)}$.

In regions $R_3$ and $R_4$ we use the split
$S_1 := \{(x_1,x_2) \in \R^2 : 0 \le x_2 \le 1\}$ and in
regions $R_5$ and $R_6$ we use the split
$S_2 := \{(x_1,x_2) \in \R^2 : 0 \le x_1 \le 1\}$.
In regions $R_1$ and $R_2$ we choose either $S_1$ or $S_2$.
We will use the one whose normal vector attains $w$,
i.e.~we choose $S_1$ if $a_2 \le 2$ and $S_2$ if $a_2 > 2$.
Note that our choice of the single split which is taken
instead of the split closure is such that it attains minimal
lattice width with respect to $T_2$.
Such a split has another nice property: among all splits
that contain $f$ in their interior we select one which
covers most of the area of the triangle.

\subsection{Regions $R_3$ and $R_4$} \label{t2.reg3reg4}

Let $f \in R_3 \cup R_4$.
In order to compute $\bar{t}(T_2, S_1)$ we have to solve
the optimization problem
\begin{equation} \label{t2.r3r4}
  \begin{split}
    \min\ & s_1 + s_2 + s_3 \\ \mbox{\quad s.t. } &
    \psi^{S_1}(r^1)s_1 + \psi^{S_1}(r^2) s_2 +
    \psi^{S_1}(r^3)s_3  \ge 1, \\
    & s_j \ge 0 \quad \mbox{for } j = 1, 2, 3.
  \end{split}
\end{equation}
Since we use $S_1$ we have $\pi = (0,1)$ and since $f \in
R_3 \cup R_4$ it holds $0 < \pi^Tf < 1$.
Thus, using \eqref{coefficients}, we obtain $\psi^{S_1}(r^1)
= \psi^{S_1}(r^2) = 1$ and $\psi^{S_1}(r^3) =
\frac{a_2-f_2}{1-f_2}$.
Therefore, an optimal solution to \eqref{t2.r3r4} is $s^\ast
= (0,0,\frac{1-f_2}{a_2-f_2})$ with optimal objective value
$\frac{1-f_2}{a_2-f_2}$.
It follows $\bar{t}(T_2, S_1) = \frac{a_2-f_2}{1-f_2}$
for $f \in R_3 \cup R_4$.

\subsection{Regions $R_5$ and $R_6$} \label{t2.reg5reg6}

Let $f \in R_5$.
To compute $\bar{t}(T_2, S_2)$ in this case we solve the
optimization problem $\min\ \sum_{j=1}^{3}{s_j}$
such that $\sum_{j=1}^{3}{\psi^{S_2}(r^j)s_j} \ge 1,\
s_j \ge 0$ for $j = 1,2,3$.
Applying \eqref{coefficients} we obtain $\psi^{S_2}(r^1) =
\frac{f_1 + \frac{a_1}{a_2-1}}{f_1}$, $\psi^{S_2}(r^2) =
\frac{\frac{a_2-a_1}{a_2-1}-f_1}{1-f_1}$, and
$\psi^{S_2}(r^3) = \frac{a_1-f_1}{1-f_1}$.
Hence, the optimal solution is the minimum among
$\{ (\psi^{S_2}(r^1))^{-1}$, $(\psi^{S_2}(r^2))^{-1}$,
$(\psi^{S_2}(r^3))^{-1}\}$.
Using our assumptions on the variables and parameters it
follows $(\psi^{S_2}(r^3))^{-1}$ $\ge 1$ and
$(\psi^{S_2}(r^i))^{-1} \le 1$ for $i \in \{1,2\}$.
One easily verifies that  $(\psi^{S_2}(r^1))^{-1} \le
(\psi^{S_2}(r^2))^{-1} \Leftrightarrow a_1 \ge f_1$ which is
satisfied in $R_5$ by assumption.
Thus, $\bar{t}(T_2, S_2) = \frac{f_1 +
\frac{a_1}{a_2-1}}{f_1}$ for $f \in R_5$.

Let $f \in R_6$.
By symmetry, i.e.~$a_1 \rightarrow 1-a_1$ and $f_1
\rightarrow 1-f_1$, we obtain $\bar{t}(T_2, S_2) =
\frac{\frac{a_2-a_1}{a_2-1}-f_1}{1-f_1}$ for $f \in R_6$.

\subsection{Regions $R_1$ and $R_2$} \label{t2.reg1reg2}

First assume $a_2 \le 2$ and use the split $S_1$.
Let $f \in R_1 \cup R_2$.
The associated optimization problem is $\min\
\sum_{j=1}^{3}{s_j}$
s.t. $\sum_{j=1}^{3}{\psi^{S_1}(r^j)s_j} \ge 1,\ s_j \ge 0$
for $j = 1,2,3$ with $\psi^{S_1}(r^1) = \psi^{S_1}(r^2) = 1$
and $\psi^{S_1}(r^3) = \frac{a_2-f_2}{1-f_2}$.
Hence, $\bar{t}(T_2, S_1) =
\frac{a_2-f_2}{1-f_2}$ for $f \in R_1 \cup R_2$ and $a_2 \le
2$.

Now assume that $a_2 > 2$.
We use the split $S_2$.
Let $f \in R_1$.
The solution of the optimization problem $\min\
\sum_{j=1}^{3}{s_j}$
s.t. $\sum_{j=1}^{3}{\psi^{S_2}(r^j)s_j} \ge 1,\ s_j \ge 0$ 
for $j = 1,2,3$ with $\psi^{S_2}(r^1) = \frac{f_1 +
\frac{a_1}{a_2-1}}{f_1}$, $\psi^{S_2}(r^2) =
\frac{\frac{a_2-a_1}{a_2-1}-f_1}{1-f_1}$, and
$\psi^{S_2}(r^3) = \frac{a_1-f_1}{1-f_1}$ is the minimum of
the set $\{ (\psi^{S_2}(r^1))^{-1},
(\psi^{S_2}(r^2))^{-1}\}$ as $(\psi^{S_2}(r^3))^{-1} \ge 1$.
It is easy to verify that $(\psi^{S_2}(r^1))^{-1} \le
(\psi^{S_2}(r^2))^{-1} \Leftrightarrow a_1 \ge f_1$.
Therefore, we have $\bar{t}(T_2, S_2) = \frac{f_1 +
\frac{a_1}{a_2-1}}{f_1}$ for $f \in R_1$ and $a_2 > 2$.

Finally, assume $f \in R_2$.
By symmetry, we obtain $\bar{t}(T_2, S_2) =
\frac{\frac{a_2-a_1}{a_2-1}-f_1}{1-f_1}$ for $f \in R_2$ and
$a_2 > 2$.

The following table summarizes the function $\bar{t}(T_2,
S_i)$, for the corresponding $i \in \{1,2\}$.

\begin{center}
$
\begin{array}{c|c|l}
  \textup{$\bar{t}(T_2, S_i)$} &
  \textup{$\bar{t}(T_2, S_i)$} &
  \textup{Location of $f$} \\
  \textup{for $a_2 \le 2$} & \textup{for $a_2 > 2$} &
  \\ \hline
  & & \\[-3mm]
  \frac{a_2-f_2}{1-f_2} &
  \frac{f_1 + \frac{a_1}{a_2-1}}{f_1} &
  f \in R_1 \\
  \frac{a_2-f_2}{1-f_2} &
  \frac{\frac{a_2-a_1}{a_2-1}-f_1}{1-f_1} &
  f \in R_2 \\
  \frac{a_2-f_2}{1-f_2} &
  \frac{a_2-f_2}{1-f_2} &
  f \in R_3 \\
  \frac{a_2-f_2}{1-f_2} & 
  \frac{a_2-f_2}{1-f_2} &
  f \in R_4 \\
  \frac{f_1 + \frac{a_1}{a_2-1}}{f_1} &
  \frac{f_1 + \frac{a_1}{a_2-1}}{f_1} &
  f \in R_5 \\
  \frac{\frac{a_2-a_1}{a_2-1}-f_1}{1-f_1} &
  \frac{\frac{a_2-a_1}{a_2-1}-f_1}{1-f_1} &
  f \in R_6
\end{array}
$
\end{center}

\subsection{Approximation for $P^{T_2}(z)$} \label{t2.integrals}

We compute the integrals
$\int_{f \in R_j}{\boldsymbol{1}\{\bar{t}(T_2,S_i) \le
z\}df}$ for $j \in \{1, \dots, 6\}$ and the corresponding
split $S_1$ or $S_2$ which we used above.
For simplicity,
$\int_{f \in R_j} := \int_{f \in
  R_j}{\boldsymbol{1}\{\bar{t} \le z\}df}$
for $j = 1, \dots, 6$.

\begin{figure}[ht]
  \unitlength=1mm
  \begin{center}
    \subfigure[\label{t2ale2} $R_1$ and $R_2$ for $a_2 \le
    2$]{
      \, \, \,
      \includegraphics[height=15\unitlength]{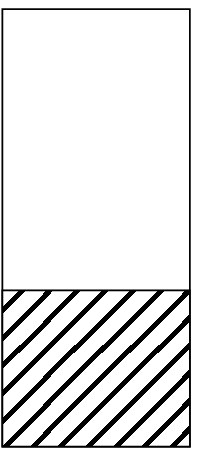}
      \, \, \,
    }
    \qquad
    \subfigure[\label{t2ag2r1} $R_1$ for $a_2 > 2$]{
      \;
      \includegraphics[height=15\unitlength]{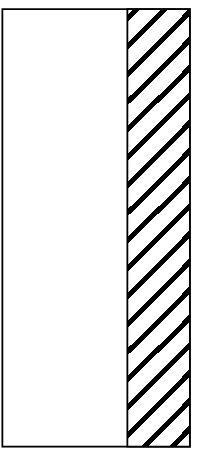}
      \;
    }
    \qquad
    \subfigure[\label{t2ag2r2} $R_2$ for $a_2 > 2$]{
      \;
      \includegraphics[height=15\unitlength]{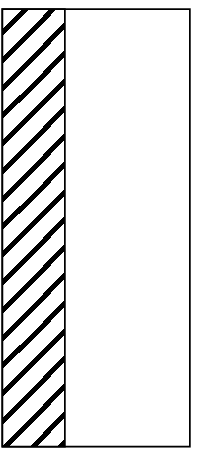}
      \;
    }
    \qquad
    \subfigure[\label{t2r3} $R_3$]{
      \includegraphics[width=15\unitlength]{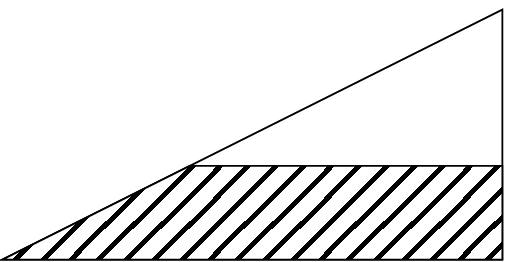}
    }
    \qquad
    \subfigure[\label{t2r4} $R_4$]{
      \includegraphics[width=15\unitlength]{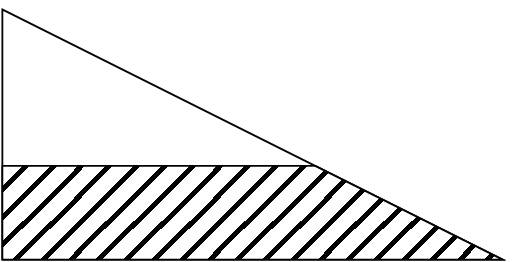}
    } 
    \qquad 
    \subfigure[\label{t2r5} $R_5$]{
      \includegraphics[height=15\unitlength]{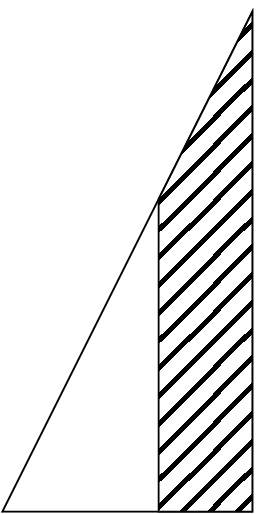}
    }
    \qquad
    \subfigure[\label{t2r6} $R_6$]{
      \includegraphics[height=15\unitlength]{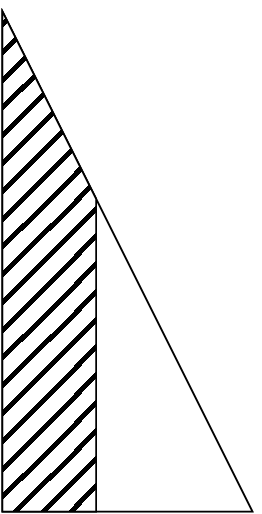}
    }
  \end{center}
  \caption{\label{t2.regions} The shaded regions satisfy
    $\boldsymbol{1}\{\bar{t} \le z\} = 1$.}
\end{figure}

Let $f \in R_3$.
Then $\boldsymbol{1}\{\bar{t} \le z\} = 1 \Leftrightarrow
\frac{a_2-f_2}{1-f_2} \le z \Leftrightarrow f_2 \le
\frac{z-a_2}{z-1}$.
Observe that $\int_{f \in R_3}
= 0$ if $z \le a_2$.
If $z > a_2$, then the area
of the set $\{f \in R_3: f_2 \le \frac{z-a_2}{z-1}\}$
is the difference of the area
of two triangles (see Fig.~\ref{t2r3}).
Direct computations yield
\begin{equation*}
  \int\limits_{f \in R_3} =
  \begin{cases}
    0 & \textup{if } z \le a_2, \\
    \frac{a_1(z-a_2)(z+a_2-2)}{2(a_2-1)(z-1)^2} &
    \textup{if } z > a_2. \\ 
  \end{cases}
\end{equation*}
Let $f \in R_5$.
Then $\boldsymbol{1}\{\bar{t} \le z\} = 1  \Leftrightarrow
\frac{f_1 + \frac{a_1}{a_2-1}}{f_1} \le z
\Leftrightarrow f_1 \ge \frac{a_1}{(a_2-1)(z-1)}$.
We have $\int_{f \in R_5}
= 0 \Leftrightarrow \frac{a_1}{(a_2-1)(z-1)} \ge a_1
\Leftrightarrow z \le \frac{a_2}{a_2-1}$.
If $z > \frac{a_2}{a_2-1}$, then again we compute the
difference of the area of two triangles (see
Fig.~\ref{t2r5}) and infer
\begin{equation*}
  \int\limits_{f \in R_5} =
  \begin{cases}
    0 & \textup{if } z \le \frac{a_2}{a_2-1}, \\
    \frac{a_1(a_2-1)}{2} g_3 & 
    \textup{if } z > \frac{a_2}{a_2-1}, \\ 
  \end{cases}
\end{equation*}
where $g_3 = 1-\frac{1}{(a_2-1)^2(z-1)^2}$.

By symmetry, the integrals for $f \in R_4$ and $f \in R_6$
are obtained by replacing $a_1$ with $1-a_1$ in the formulas
for $R_3$ and $R_5$ (see Fig.~\ref{t2r4} and \ref{t2r6}).
Thus,
\begin{align*}
  \int\limits_{f \in R_4} &=
  \begin{cases}
    0 & \textup{if } z \le a_2, \\
    \frac{(1-a_1)(z-a_2)(z+a_2-2)}{2(a_2-1)(z-1)^2} &
    \textup{if } z > a_2, \\ 
  \end{cases} \\
  \int\limits_{f \in R_6} &=
  \begin{cases}
    0 & \textup{if } z \le \frac{a_2}{a_2-1}, \\
    \frac{(1-a_1)(a_2-1)}{2} g_3 & 
    \textup{if } z > \frac{a_2}{a_2-1}. \\ 
  \end{cases}
\end{align*}
In order to compute
$\int_{f \in R_i}$ for
$i = 1,2$ we
distinguish the cases $a_2 \le 2$ and $a_2 > 2$.
First, let us assume that $a_2 \le 2$.
Computations in a similar manner as above (area of the
shaded  quadrilateral in Fig.~\ref{t2ale2}) lead to
\begin{align*}
  \int\limits_{f \in R_1} &=
  \begin{cases}
    0 & \textup{if } z \le a_2, \\
    \frac{a_1(z-a_2)}{z-1} & \textup{if } z > a_2, \\ 
  \end{cases} \\
  \int\limits_{f \in R_2} &=
  \begin{cases}
    0 & \textup{if } z \le a_2, \\
    \frac{(1-a_1)(z-a_2)}{z-1} & \textup{if } z > a_2. \\ 
  \end{cases}
\end{align*}
Now let us assume that $a_2 > 2$.
Computing the corresponding area (see Fig.~\ref{t2ag2r1} and
\ref{t2ag2r2}) yields
\begin{align*}
  \int\limits_{f \in R_1} &=
  \begin{cases}
    0 & \textup{if } z \le \frac{a_2}{a_2-1}, \\
    a_1 g_4 &
    \textup{if } z > \frac{a_2}{a_2-1}, \\ 
  \end{cases} \\
  \int\limits_{f \in R_2} &=
  \begin{cases}
    0 & \textup{if } z \le \frac{a_2}{a_2-1}, \\
    (1-a_1) g_4 &
    \textup{if } z > \frac{a_2}{a_2-1}, \\ 
  \end{cases}
\end{align*}
where $g_4 = 1 - \frac{1}{(a_2-1)(z-1)}$.

We can aggregate the formulas for regions $R_1$ and $R_2$,
$R_3$ and $R_4$, and $R_5$ and $R_6$ in order to
eliminate the parameter $a_1$.
It follows
\begin{align*}
  \int\limits_{f \in R_1} +
  \int\limits_{f \in R_2} &=
  \begin{cases}
    0 & \textup{if } z \le a_2 \textup{ and } a_2 \le 2, \\
    \frac{z-a_2}{z-1} & \textup{if } z > a_2 \textup{ and }
    a_2 \le 2, \\ 
  \end{cases} \\
  \int\limits_{f \in R_1} +
  \int\limits_{f \in R_2} &=
  \begin{cases}
    0 & \textup{if } z \le \frac{a_2}{a_2-1} \textup{ and }
    a_2 > 2, \\
    1 - \frac{1}{(a_2-1)(z-1)}
    & \textup{if } z > \frac{a_2}{a_2-1} \textup{ and }
    a_2 > 2, \\ 
  \end{cases} \\
  \int\limits_{f \in R_3} +
  \int\limits_{f \in R_4} &=
  \begin{cases}
    0 & \textup{if } z \le a_2, \\
    \frac{(z-a_2)(z+a_2-2)}{2(a_2-1)(z-1)^2}
    & \textup{if } z > a_2, \\ 
  \end{cases} \\
  \int\limits_{f \in R_5} +
  \int\limits_{f \in R_6} &=
  \begin{cases}
    0 & \textup{if } z \le \frac{a_2}{a_2-1}, \\
    \frac{a_2-1}{2} \left(1 - \frac{1}{(a_2-1)^2(z-1)^2}\right) &
    \textup{if } z > \frac{a_2}{a_2-1}. \\ 
  \end{cases}
\end{align*}
We are now ready to state our results in terms of the lattice
width.
Reinterpreting the formulas above we obtain for $a_2 \le 2$
(i.e.~$w = a_2$) the integrals
\begin{align*}
  \int\limits_{f \in R_1} +
  \int\limits_{f \in R_2} &=
  \begin{cases}
    0 & \textup{if } z \le w, \\
    \frac{z-w}{z-1} & \textup{if } z > w, \\ 
  \end{cases} \\
  \int\limits_{f \in R_3} +
  \int\limits_{f \in R_4} &=
  \begin{cases}
    0 & \textup{if } z \le w, \\
    \frac{(z-w)(z+w-2)}{2(w-1)(z-1)^2}
    & \textup{if } z > w, \\ 
  \end{cases} \\
  \int\limits_{f \in R_5} +
  \int\limits_{f \in R_6} &=
  \begin{cases}
    0 & \textup{if } z \le \frac{w}{w-1}, \\
    \frac{(w-1)^2(z-1)^2-1}{2(w-1)(z-1)^2} &
    \textup{if } z > \frac{w}{w-1}. \\ 
  \end{cases}
\end{align*}
For $a_2 > 2$ (i.e.~$w = \frac{a_2}{a_2-1}$) we have
\begin{align*}
  \int\limits_{f \in R_1} +
  \int\limits_{f \in R_2} &=
  \begin{cases}
    0 & \textup{if } z \le w, \\
    \frac{z-w}{z-1} & \textup{if } z > w, \\ 
  \end{cases}
\end{align*}


\begin{align*}
  \int\limits_{f \in R_3} +
  \int\limits_{f \in R_4} &=
  \begin{cases}
    0 & \textup{if } z \le \frac{w}{w-1}, \\
     \frac{(w-1)^2(z-1)^2-1}{2(w-1)(z-1)^2}
    & \textup{if } z > \frac{w}{w-1}, \\ 
  \end{cases} \\
  \int\limits_{f \in R_5} +
  \int\limits_{f \in R_6} &=
  \begin{cases}
    0 & \textup{if } z \le w, \\
    \frac{(z-w)(z+w-2)}{2(w-1)(z-1)^2} &
    \textup{if } z > w. \\ 
  \end{cases}
\end{align*}
Note that $w \le \frac{w}{w-1}$ as $w \le 2$ (see
\cite{AverkovWagner} for a proof that any triangle of type
2 satisfies $1 < w \le 2$).
Hence, together with \eqref{T2.approx} we arrive at the
following formula.
\begin{equation} \label{lwidth.formula}
  P^{T_2}(z) \ge
  \begin{cases}
    0 & \textup{if } 1 < z \le w, \\[1.8mm]
    g_1 &
    \textup{if } w < z \le \frac{w}{w-1}, \\[1.8mm]
    g_1 + g_2 &
    \textup{if } \frac{w}{w-1} < z < + \infty, \\
  \end{cases}
\end{equation}
with $g_1 = \frac{(z-w)(2wz-w-z)}{w^2(z-1)^2}$ and
$g_2 = \frac{(w-1)^2(z-1)^2-1}{w^2(z-1)^2}$.
Let $z > 1$ be given.
We want to show that the probability to improve upon the
split closure by a value of more than $z$ when adding a type
2 triangle becomes the smaller the closer the type 2
triangle is to a split, i.e.~the closer $w$ is to $1$.
For $w$ being sufficiently close to $1$ we have $w < z \le
\frac{w}{w-1}$.
Just substitute $1$ for $w$ in \eqref{lwidth.formula} to
infer that $P^{T_2}(z)$ converges to $1$ for any $z > 1$.
Therefore, $1 - P^{T_2}(z)$ tends to $0$.
In other words, the chance that adding a flat type 2
triangle improves the split closure by a value of more than
$z$ with respect to our strength measure tends to $0$.
In terms of the vertex $(a_1,a_2)$ this happens if $a_2$
converges either to $1$ (i.e.~$T_2$ converges to $S_1$) or
infinity (i.e.~$T_2$ converges to $S_2$).

\begin{figure}[ht]
  \unitlength=1mm
  \begin{center}
    \subfigure[\label{z=2} An upper bound for $1 -
    P^{T_2}(2)$ in terms of $w$.]{
      \includegraphics[height=30\unitlength]{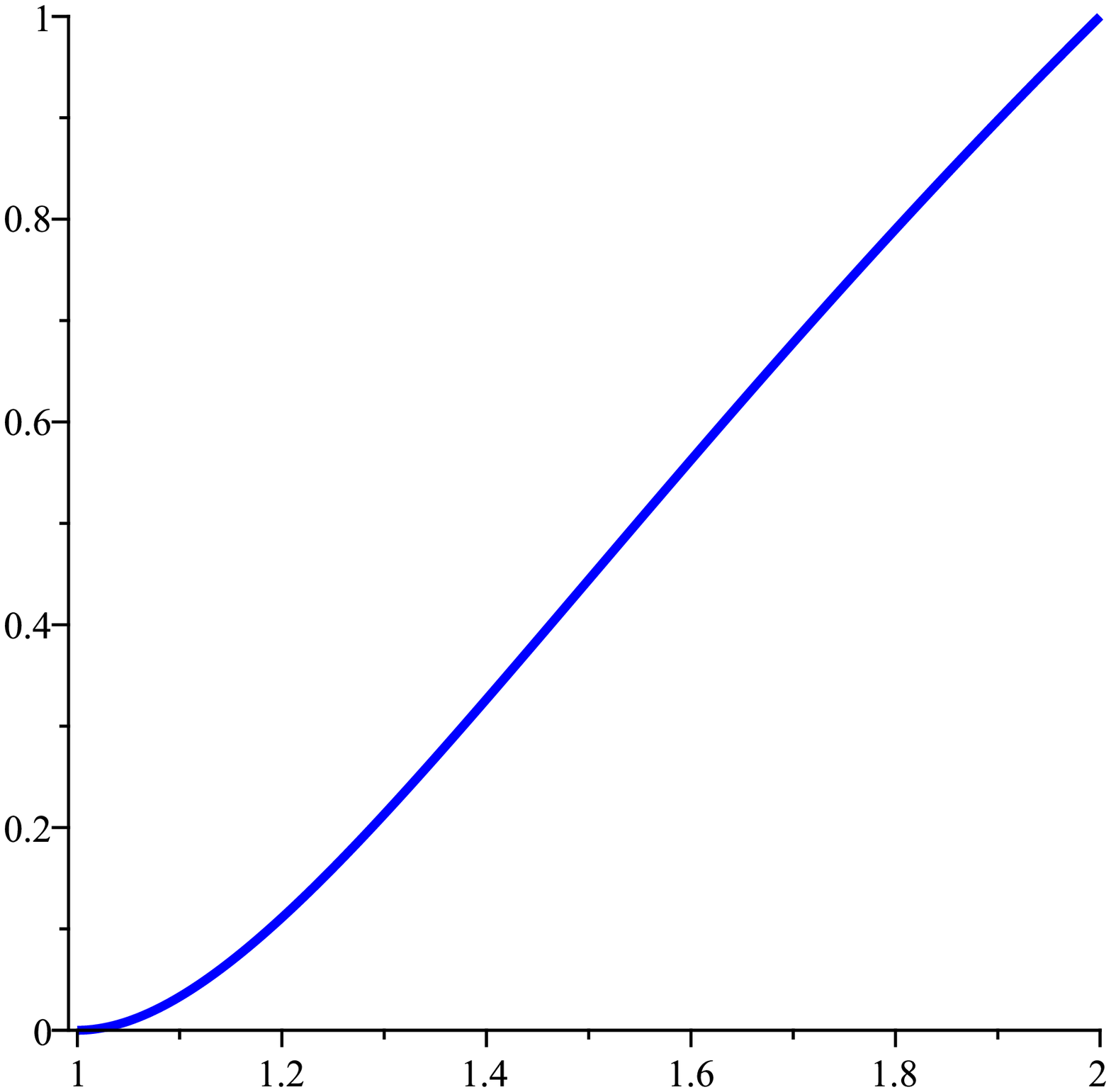}
    }
    \qquad
    \subfigure[\label{z=3/2} A lower bound for
    $P^{T_2}(\frac{3}{2})$ in terms of $w$.]{
      \includegraphics[height=30\unitlength]{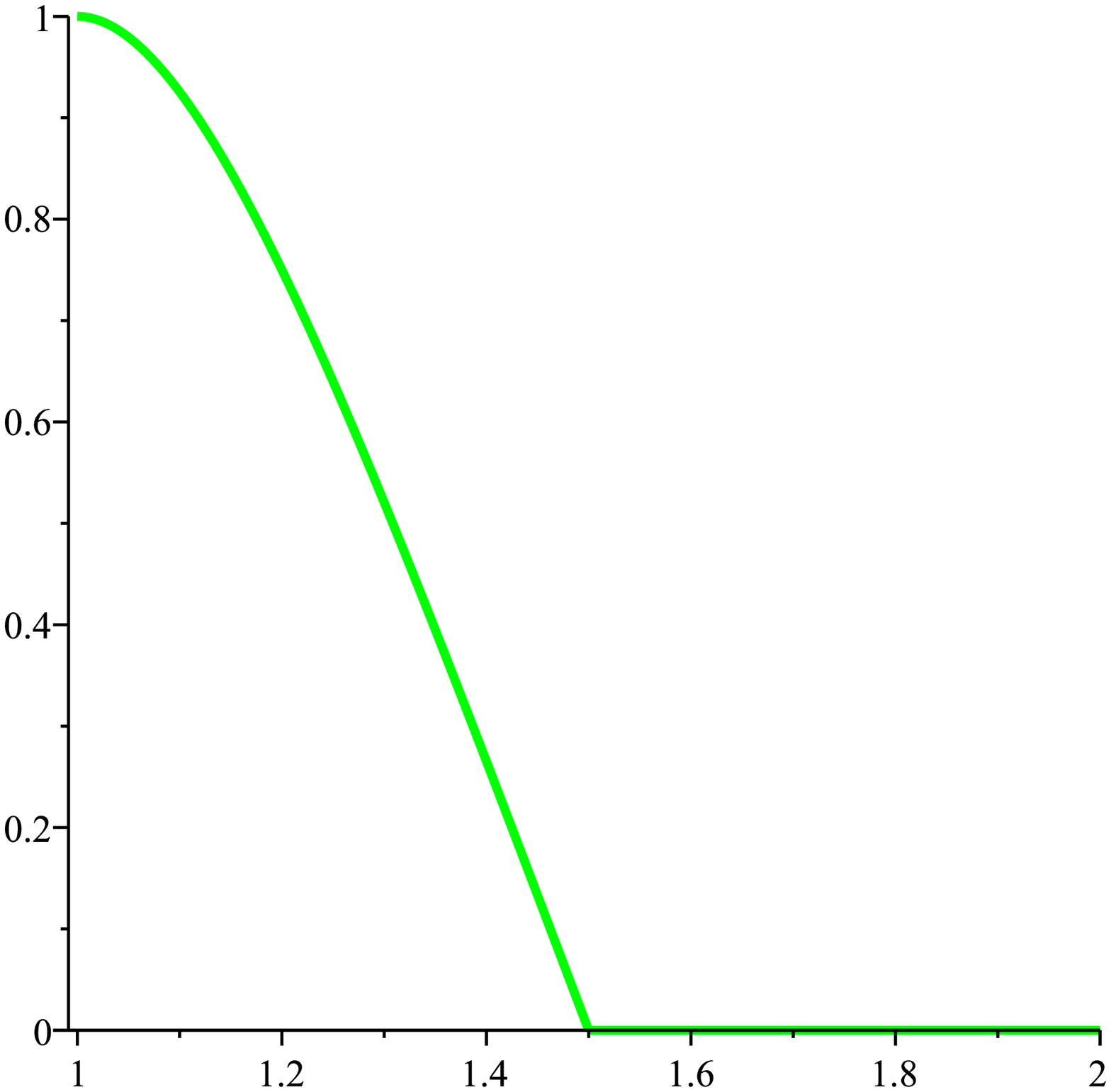}
    }
  \end{center}
  \caption{Bounds for $1 - P^{T_2}(2)$ and $P^{T_2}(\frac{3}{2})$.}
\end{figure}

Two special cases are of interest: $z = 2$ and $z =
\frac{3}{2}$.
Plugging in $z = 2$ in \eqref{lwidth.formula} yields
$1 - P^{T_2}(2) \le \frac{4(w-1)^2}{w^2}$ for any $1 < w \le
2$ (see Fig.~\ref{z=2}) and plugging in $z = \frac{3}{2}$ in
\eqref{lwidth.formula} leads to
\begin{equation*}
  P^{T_2}\left(\frac{3}{2}\right) \ge
  \begin{cases}
    0 & \textup{if } \frac{3}{2} \le w \le 2, \\
    \frac{(3-2w)(4w-3)}{w^2} &
    \textup{if } 1 < w < \frac{3}{2} \\ 
  \end{cases}
\end{equation*}
(see Fig.~\ref{z=3/2}).
$1 - P^{T_2}(2)$ is interpretable as the probability that
the inequality corresponding to $T_2$ closes more gap to
$R_f^3$ than the inequality corresponding to any type 1
triangle.
On the other hand $P^{T_2}(\frac{3}{2})$ can be seen as the
probability that adding $T_2$ to the split closure is
inferior to adding any type 1 triangle $T_1$ with rays going
through the corners of $T_1$.

The interpretation is based on the fact that any type 1
triangle (assuming corner rays) has a strength between
$\frac{3}{2}$ and $2$.
On the one hand type 1 triangles are inferior to flat type 2
triangles by comparing the worst case strength (see
\cite{BasuBonamiCornuejolsMargot09}).
On the other hand a type 1 triangle guarantees a strength
of $\frac{3}{2}$ and is therefore superior to a flat type 2
triangle on average.
For instance, let $T_2$ be a type 2 triangle with $w(T_2) =
1{.}1$.
Then $1 - P^{T_2}(2) < 3{.}4 \%$ and
$P^{T_2}(\frac{3}{2}) > 92{.}5 \%$.
Thus, with a probability of less than $3{.}4 \%$, $T_2$ is
better than any type 1 triangle $T_1$, and worse with a
probability of more than $92{.}5 \%$ if corner rays are
assumed for $T_1$.
We point out that these probabilities are quite close to
$0$ and $1$, respectively, even though we used just one
split instead of the entire split closure.


\section{Quadrilaterals} \label{sec.quad}

\begin{figure}[ht]
        \centerline{\includegraphics[scale = 1.0]{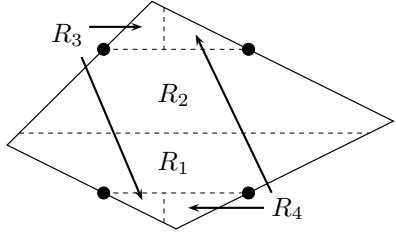}}
        \caption{Decomposition of a quadrilateral.}
        \label{quad}
\end{figure}

The vertices of the quadrilateral $Q$ are denoted $a =
(a_1,a_2)$, $b = (b_1,b_2)$, $c = (c_1,c_2)$, and $d =
(d_1,d_2)$.
By an affine unimodular transformation, we assume that the
point $(0,0)$
(resp.~$(1,0)$, $(0,1)$, $(1,1)$) is in the relative
interior of the facet with vertices $b$ and $c$ (resp.~$b$
and $d$, $a$ and $c$, $a$ and $d$).
We further assume $0 < a_1 \le b_1 < 1$, $1 < a_2$, $b_2 <
0$, and $-b_2 \le a_2 - 1$.
The parameters $a_1$, $a_2$, $b_1$, and $b_2$ are assumed to
be arbitrary but fixed.
This implies $c_1 = -\frac{a_1b_1}{(a_2-1)b_1 - a_1b_2}$,
$c_2 = \frac{c_1b_2}{b_1}$, $d_1 =
\frac{(a_2-a_1)(1-b_1)-(1-a_1)b_2}{(a_2-1)(1-b_1)-(1-a_1)b_2}$,
and $d_2 = \frac{(d_1-1)b_2}{b_1-1}$.
One easily verifies $c_1 < 0$, $0< c_2 < 1$, $1< d_1$, $0 <
d_2 < 1$, $c_2 \le d_2$, and $\area(Q) = \frac{1}{2}(a_2 -
b_2 + d_1 - c_1)$.
Under these assumptions it is known that $w := w(Q) =
\min\{a_2-b_2,d_1-c_1\}$ (see \cite{AverkovWagner}).
Thus, without loss of generality we assume $w = a_2 - b_2$.
We decompose $Q$ into four regions:
$R_1 := \intt(Q \cap \{(x_1,x_2) \in \R^2 :
0 \le x_2 \le \frac{-b_2}{w-1}\})$, 
$R_2 := \intt(Q \cap \{(x_1,x_2) \in \R^2 :
\frac{-b_2}{w-1} \le x_2 \le 1\})$,
$R_3 := \intt((Q \setminus \{R_1 \cup R_2\}) \cap
\{(x_1,x_2) \in \R^2 :
0 \le x_1 \le \theta\})$, 
$R_4 := \intt((Q \setminus \{R_1 \cup R_2\}) \cap
\{(x_1,x_2) \in \R^2 :
\theta \le x_1 \le 1\})$, where $\theta =
\frac{a_1b_1(a_2-1)(1-b_1) -
  b_1a_1(1-a_1)b_2}{b_1(a_2-1)(1-b_1) - a_1(1-a_1)b_2}$ (see
Fig.~\ref{quad}).
It is tedious but easy to verify that
$c_2 \le \frac{-b_2}{w-1} \le d_2$ and 
$a_1 \le \theta \le b_1$.
As in the case of type 2 triangles, instead of taking the
split closure, we use a single well-chosen split inequality.

In regions $R_1$ and $R_2$ we use the split $S_1 :=
\{(x_1,x_2) \in \R^2 : 0 \le x_2 \le 1\}$ and in regions
$R_3$ and $R_4$ the split $S_2 := \{(x_1,x_2) \in \R^2 : 0
\le x_1 \le 1\}$.
Since the computations of $\bar{t}(Q, S_1)$ for $f \in
R_1 \cup R_2$ and $\bar{t}(Q, S_2)$ for $f \in R_3 \cup
R_4$ are straightforward (see Sections \ref{t2.reg3reg4} --
\ref{t2.reg1reg2} for an illustration on how to do that for
type 2 triangles) we only state the results.

\begin{center}
$
\begin{array}{c|c}
  \textup{$\bar{t}(Q, S_i)$} &
  \textup{Location of $f$} \\ \hline
  & \\[-3mm]
  \frac{f_2-b_2}{f_2} &
  f \in R_1 \\
  \frac{a_2-f_2}{1-f_2} &
  f \in R_2 \\
  \frac{f_1-c_1}{f_1} &
  f \in R_3 \\
  \frac{d_1-f_1}{1-f_1} &
  f \in R_4 \\
\end{array}
$
\end{center}

As in Section \ref{t2.integrals}, we now compute the
integrals in order to obtain a lower bound for $P^{Q}(z)$.
For simplicity, let
$\int_{f \in R_j} := \int_{f \in
  R_j}{\boldsymbol{1}\{\bar{t} \le z\}df}$
for $j = 1, \dots, 4$.

\subsection{Regions $R_1$ and $R_2$}

Let $f \in R_1$.
Then $\boldsymbol{1}\{\bar{t} \le z\} = 1
\Leftrightarrow \frac{f_2-b_2}{f_2} \le z \Leftrightarrow
f_2 \ge \frac{-b_2}{z-1}$.
Observe that $\int_{f \in R_1}
= 0$ if $\frac{-b_2}{w-1} \le \frac{-b_2}{z-1} \Leftrightarrow
z \le w$.
If $z > w$, then we distinguish into $\frac{-b_2}{z-1} \ge
c_2$ and $\frac{-b_2}{z-1} \le c_2$.
If $\frac{-b_2}{z-1} \ge c_2$, then the area to compute is a
trapezoid with area $A_1$;
if $\frac{-b_2}{z-1} \le c_2$ the area to compute is
composed of two trapezoids with aggregate area $A_2$,
where
\begin{align*}
  A_1 = \frac{1}{2}& \left( \frac{-b_2}{w-1} -
    \frac{-b_2}{z-1} \right)\left( \frac{w-(b_1-a_1)}{w-1} +
    \frac{z-b_1}{z-1} + \frac{a_1(z-1+b_2)}{(a_2-1)(z-1)}
    \right), \\
  A_2 = \frac{1}{2}& \left( \frac{-b_2}{w-1} - c_2
    \right)\left( \frac{w-(b_1-a_1)}{w-1} + \frac{a_1(b_2-1)
    - (a_2-1)b_1}{a_1b_2-(a_2-1)b_1} \right) \\
     &+ \frac{1}{2} \left( c_2 - \frac{-b_2}{z-1}
    \right)\left( \frac{z}{z-1} + \frac{a_1(b_2-1)
    - (a_2-1)b_1}{a_1b_2-(a_2-1)b_1} \right).
\end{align*}
We obtain 
$$
\int\limits_{f \in R_1} =
  \begin{cases}
    0 & \textup{if } 1 < z < w, \\
    A_1 & \textup{if } w \le z \le \frac{c_2-b_2}{c_2}, \\
    A_2 & \textup{if } \frac{c_2-b_2}{c_2} < z < + \infty. \\
  \end{cases}
$$

Let $f \in R_2$.
Then $\boldsymbol{1}\{\bar{t} \le z\} = 1
\Leftrightarrow \frac{a_2-f_2}{1-f_2} \le z \Leftrightarrow
f_2 \le \frac{z-a_2}{z-1}$.
Thus, $\int_{f \in R_2}
= 0$ if $\frac{z-a_2}{z-1} \le \frac{-b_2}{w-1} \Leftrightarrow
z \le w$.
Otherwise, we distinguish into $\frac{z-a_2}{z-1} \le
d_2$ and $\frac{z-a_2}{z-1} \ge d_2$.
In the first case the area to compute is a trapezoid with
area $A_3$; 
in the second case the area is composed of two trapezoids
with aggregate area $A_4$, where
\begin{align*}
  A_3 = \frac{1}{2}& \left( \frac{z-a_2}{z-1} -
    \frac{-b_2}{w-1} \right)\left( \frac{w-(b_1-a_1)}{w-1} +
    \frac{z-1+a_1}{z-1} + \frac{(z-a_2)(b_1-1)}{b_2(z-1)}
    \right), \\
  A_4 = \frac{1}{2}& \left( \frac{z-a_2}{z-1} - d_2
    \right)\left( \frac{z}{z-1} +
      \frac{a_2(1-b_1)-(1-a_1)b_2}{(a_2-1)(1-b_1)-(1-a_1)b_2}
    \right) \\
     &+ \frac{1}{2} \left( d_2 - \frac{-b_2}{w-1}
    \right)\left( \frac{w-(b_1-a_1)}{w-1} +
      \frac{a_2(1-b_1)-(1-a_1)b_2}{(a_2-1)(1-b_1)-(1-a_1)b_2} 
    \right).
\end{align*}
We infer
$$
\int\limits_{f \in R_2} =
  \begin{cases}
    0 & \textup{if } 1 < z < w, \\
    A_3 & \textup{if } w \le z \le \frac{a_2-d_2}{1-d_2}, \\
    A_4 & \textup{if } \frac{a_2-d_2}{1-d_2} < z < + \infty. \\
  \end{cases}
$$

\subsection{Regions $R_3$ and $R_4$}

Applying the same procedure as in the previous section
leads to
\begin{align*}
  A_5 = \frac{1}{2}& \left( \frac{-c_1}{d_1-c_1-1} -
    \frac{-c_1}{z-1} \right)\left(
    \frac{(a_2-1)(d_1-1)}{(1-a_1)(d_1-c_1-1)} +
    \frac{(a_2-1)(z-1+c_1)}{(1-a_1)(z-1)} +
    \frac{c_2}{d_1-c_1-1} + \frac{c_2}{z-1}
    \right), \\
  A_6 = \frac{1}{2}& \left( \frac{-c_1}{d_1-c_1-1} - a_1
    \right)\left( \frac{(a_2-1)(2-a_1)}{1-a_1} -
      \frac{a_1b_2}{b_1} +
      \frac{c_1(a_2-1)+c_2(1-a_1)}{(1-a_1)(d_1-c_1-1)}
    \right) \\
     &+ \frac{1}{2} \left( a_1 - \frac{-c_1}{z-1}
    \right)\left( a_2 - 1 - \frac{a_1b_2}{b_1} +
      \frac{a_1c_2 - c_1(a_2-1)}{a_1(z-1)} 
    \right), \\
  A_7 = \frac{1}{2}& \cdot
    \frac{(d_1-1)(z-d_1+c_1)}{(z-1)(d_1-c_1-1)} \cdot
    \left( \frac{c_2(1-a_1)+(a_2-1)(d_1-1)}{(1-a_1)(d_1-c_1-1)} + 
    \frac{(a_2-1)(d_1-1)}{(1-a_1)(z-1)} -
    \frac{b_2(z-d_1)}{b_1(z-1)}
    \right), \\
  A_8 = \frac{1}{2}& \left( b_1 - \frac{-c_1}{d_1-c_1-1}
    \right)\left( \frac{c_2(1-a_1) +
    c_1(a_2-1)}{(1-a_1)(d_1-c_1-1)} + 
    \frac{(a_2-1)(2-b_1)-b_2(1-a_1)}{1-a_1}
    \right) \\
     &+ \frac{1}{2} \left( \frac{z-d_1}{z-1} - b_1
    \right)\left( \frac{(a_2-1)(z-d_1)}{(a_1-1)(z-1)} -
      \frac{b_2(d_1-1)}{(1-b_1)(z-1)} +
      \frac{(a_2-1)(2-b_1)-b_2(1-a_1)}{1-a_1}
    \right),
\end{align*}
and finally 
$$
\int\limits_{f \in R_3} =
  \begin{cases}
    0 & \textup{if } 1 < z < d_1-c_1, \\
    A_5 & \textup{if } d_1-c_1 \le z \le \frac{a_1-c_1}{a_1}, \\
    A_6 & \textup{if } \frac{a_1-c_1}{a_1} < z < + \infty, \\
  \end{cases}
$$

$$
\int\limits_{f \in R_4} =
  \begin{cases}
    0 & \textup{if } 1 < z < d_1-c_1, \\
    A_7 & \textup{if } d_1-c_1 \le z \le \frac{d_1-b_1}{1-b_1}, \\
    A_8 & \textup{if } \frac{d_1-b_1}{1-b_1} < z < + \infty. \\
  \end{cases}
$$

\subsection{Approximation for $P^Q(z)$}

Adding the integrals and dividing by the area of the
quadrilateral gives the lower bound for $P^Q(z)$ which we
wanted.
Thus,
\begin{equation} \label{quad.bound}
  P^{Q}(z) \ge \frac{1}{\area(Q)} \sum_{j=1}^{4}{
  \int\limits_{f \in R_j}{\boldsymbol{1}\{\bar{t} \le z\}df}}.
\end{equation}
We did not succeed in showing algebraically that this lower
bound tends to $1$ when $w$ converges to $1$.
However, we performed simulations supporting our
conjecture.
In these simulations we let $w$ converge to $1$ for several
values for $a_1$, $a_2$, $b_1$, and $b_2$.
Concretely, we discretized $a_1$, $a_2$, $b_1$, and $b_2$
within their ranges:
$a_1 \in \{0{.}001, 0{.}002, \dots, 0{.}999\}$,
$b_1 \in \{a_1, a_1 + 0{.}001, \dots, 0{.}999\}$,
$a_2 \in \{1{.}999, 1{.}998, \dots, 1{.}001\}$, and
$b_2 \in \{-(a_2-1), -(a_2-1) + 0{.}001, \dots, -0{.}001\}$.
In all cases we observed that the lower bound became close
to one.
This is not a proof, but gives strong indication that this
should hold in general.

We point out that under the additional assumption $a_1 =
b_1$ (implying $a_1 = \theta = b_1$ and $c_2 =
\frac{-b_2}{w-1} = d_2$) the lower bound \eqref{quad.bound}
coincides with the lower bound of a type $2$ triangle which
is stated in \eqref{lwidth.formula}.
This can be seen by simply substituting $a_1$ for $b_1$ in
the formulas for $A_1, \dots, A_8$.

\begin{corollary}
Let $Q^=$ be a quadrilateral meeting the assumptions
stated at the beginning of this section which, in addition, satisfies
$a_1 = b_1$. Let $w := w(Q^=)$.
Then 
\begin{equation}
  P^{Q^=}(z) \ge
  \begin{cases}
    0 & \textup{if } 1 < z \le w, \\[1.8mm]
    \frac{(z-w)(2wz-w-z)}{w^2(z-1)^2} &
    \textup{if } w < z \le \frac{w}{w-1}, \\[1.8mm]
    \frac{(z-w)(2wz-w-z)+(w-1)^2(z-1)^2-1}{w^2(z-1)^2} &
    \textup{if } \frac{w}{w-1} < z < + \infty. \\
  \end{cases}
\end{equation}
Moreover, for any $z > 1$, $P^{Q^=}(z)$ tends to $1$ if $w$
converges to $1$.
\end{corollary}


\section{Type 3 triangles} \label{sec.t3}

\begin{figure}[ht]
        \centerline{\includegraphics[scale = 1.0]{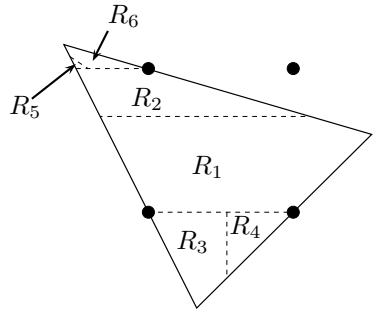}}
        \caption{Decomposition of a type 3 triangle.}
        \label{type3}
\end{figure}

By an affine unimodular transformation, we assume that the
type 3 triangle $T_3$ satisfies $T_3 \cap
\Z^2 = \{(0,0), (1,0), (0,1)\}$ and that each facet of $T_3$
contains one of these points in its relative interior.
The three vertices are denoted by $a = (a_1,a_2)$, $b =
(b_1,b_2)$, and $c = (c_1,c_2)$.
We further assume $1 < a_1$, $0 < a_2 < 1$, $0 < b_1 < 1$,
and $b_1 + b_2 < 0$.
Let $a_1$, $a_2$, and $b_1$ be arbitrary but fixed.
Thus, the other parameters are
$b_2 = -\frac{a_2(1-b_1)}{a_1-1}$,
$c_1 = \frac{a_1(a_1-1)b_1}{(a_1-1)(1-a_2)b_1-a_1a_2(1-b_1)}$, and 
$c_2 = -\frac{a_1a_2(1-b_1)}{(a_1-1)(1-a_2)b_1-a_1a_2(1-b_1)}$.
One easily verifies $b_2 < 0$, $c_1 < 0$, $1 < c_2$, $0 <
c_1 + c_2 < 1$, and $\area(T_3) =
\frac{1}{2}(a_1+a_2-b_2-c_1)$.
Under these assumptions we have $w := w(T_3) = \min\{c_2 -
b_2, a_1 - c_1, a_1 + a_2 - (b_1 + b_2)\}$ (see
\cite{AverkovWagner}).
Without loss of generality we assume $w = c_2 - b_2 \le a_1
- c_1 \le a_1 + a_2 - (b_1 + b_2)$.
During this section we consider the three splits
$S_1 := \{(x_1,x_2) \in \R^2 : 0 \le x_2 \le 1\}$,
$S_2 := \{(x_1,x_2) \in \R^2 : 0 \le x_1 \le 1\}$, and
$S_3 := \{(x_1,x_2) \in \R^2 : 0 \le x_1 + x_2 \le 1\}$.
We decompose $T_3$ into six regions:
$R_1 := \intt(T_3 \cap \{(x_1,x_2) \in \R^2 :
0 \le x_2 \le \frac{-b_2}{c_2-1-b_2}\})$, 
$R_2 := \intt(T_3 \cap \{(x_1,x_2) \in \R^2 :
\frac{-b_2}{c_2-1-b_2} \le x_2 \le 1\})$,
$R_3 := \intt((T_3 \setminus \{R_1 \cup R_2\}) \cap
\{(x_1,x_2) \in \R^2 :
0 \le x_1 \le \frac{-c_1}{a_1-1-c_1}\})$, 
$R_4 := \intt((T_3 \setminus \{R_1 \cup R_2\}) \cap
\{(x_1,x_2) \in \R^2 :
\frac{-c_1}{a_1-1-c_1} \le x_1 \le 1\})$,
$R_5 := \intt((T_3 \setminus \cup_{j=1}^4{R_j}) \cap
\{(x_1,x_2) \in \R^2 : 0 \le x_1 + x_2 \le
\frac{-(b_1+b_2)}{a_1+a_2-1-(b_1+b_2)}\})$, and
$R_6 := \intt((T_3 \setminus \cup_{j=1}^4{R_j}) \cap
\{(x_1,x_2) \in \R^2 :
\frac{-(b_1+b_2)}{a_1+a_2-1-(b_1+b_2)} \le x_1 + x_2 \le
1\})$ (see Fig.~\ref{type3}).
We point out that $R_5$ could be empty.
It is tedious but easy to verify that
$a_2 < \frac{-b_2}{c_2-1-b_2} < 1$,
$b_1 < \frac{-c_1}{a_1-1-c_1} < 1$, and
$0 < \frac{-(b_1+b_2)}{a_1+a_2-1-(b_1+b_2)} < c_1 + c_2$.

For each region, we use a single split inequality to
approximate the split closure.
In regions $R_1$ and $R_2$ we use the split $S_1$, in
regions $R_3$ and $R_4$ the split $S_2$, and in regions
$R_5$ and $R_6$ the split $S_3$.
Thus, in each region $R_j$ we choose a split which covers
$R_j$ and the convex hull of $(0,0)$, $(1,0)$, and $(0,1)$.
The following table states the values for $\bar{t}(T_3,S_i)$
for the regions $R_1$ to $R_6$. 

\begin{center}
$
\begin{array}{c|c}
  \textup{$\bar{t}(T_3, S_i)$} &
  \textup{Location of $f$} \\ \hline
  & \\[-3mm]
  \frac{f_2-b_2}{f_2} &
  f \in R_1 \\
  \frac{c_2-f_2}{1-f_2} &
  f \in R_2 \\
  \frac{f_1-c_1}{f_1} &
  f \in R_3 \\
  \frac{a_1-f_1}{1-f_1} &
  f \in R_4 \\
  \frac{f_1+f_2-(b_1+b_2)}{f_1+f_2} &
  f \in R_5 \\
  \frac{a_1+a_2-(f_1+f_2)}{1-(f_1+f_2)} &
  f \in R_6
\end{array}
$
\end{center}

The computation of the integrals $\int_{f \in
R_j}{\boldsymbol{1}\{\bar{t} \le z\}df}$ for $j = 1, 
\dots, 6$ gives no new insights as well.
We only state the results.
For simplicity, let
$\int_{f \in R_j} := \int_{f \in
  R_j}{\boldsymbol{1}\{\bar{t} \le z\}df}$
for $j = 1, \dots, 6$.
Let
\begin{align*}
  A_1 &= \frac{1}{2} \left( \frac{-b_2}{w-1} -
    \frac{-b_2}{z-1} \right)\left( \frac{b_1}{w-1} +
    \frac{b_1}{z-1} + \frac{a_1}{1-a_2} \left(
      \frac{c_2-1}{w-1} + \frac{z-1+b_2}{z-1} \right)
    \right), \\
  A_2 &= \frac{1}{2} \left( \frac{-b_2}{w-1} - a_2
    \right)\left( \frac{(1-a_2)b_1 +
        a_1(c_2-1)}{(1-a_2)(w-1)} - \frac{a_2b_1-
        a_1b_2)}{b_2} \right) \\
     &\ \ + \frac{1}{2} \left( a_2 - \frac{-b_2}{z-1}
    \right)\left( \frac{a_2b_1 - (a_1-1)b_2}{a_2(z-1)} -
      \frac{a_2b_1 - (a_1+1)b_2)}{b_2} \right), \\
  A_3 &= \frac{1}{2} \left( \frac{z-c_2}{z-1} -
    \frac{-b_2}{w-1} \right)\left( \frac{b_1}{w-1} -
    \frac{b_1(z-c_2)}{b_2(z-1)} + \frac{a_1}{1-a_2} \left(
      \frac{c_2-1}{w-1} + \frac{c_2-1}{z-1} \right)
    \right), \\
  A_4 &= \frac{1}{2} \left( \frac{-c_1}{a_1-c_1-1} -
    \frac{-c_1}{z-1} \right)\left( \frac{a_2}{a_1-c_1-1}
    + \frac{a_2(z-1+c_1)}{(a_1-1)(z-1)}
    \right), \\
  A_5 &= \frac{1}{2} \cdot \frac{a_2}{b_1(a_1-1)} \cdot
    \left( \frac{b_1(z-1)+c_1}{z-1} \right)^2, \\
  A_6 &= \frac{1}{2} \left( \frac{z-a_1}{z-1} -
    \frac{-c_1}{a_1-c_1-1} \right)\left( \frac{a_2}{a_1-c_1-1}
    + \frac{a_2}{z-1} \right), \\
  A_7 &= \frac{1}{2} \cdot
    \left( \frac{-b_2}{a_1+a_2-(b_1+b_2)-1} -1 \right)^2, \\
  A_8 &= \frac{1}{2} \left( \frac{b_1}{a_1+a_2-(b_1+b_2)-1}
    + \frac{b_1}{b_2} \right)
    \left( \frac{-b_2}{a_1+a_2-(b_1+b_2)-1} -1 \right), \\
  A_9 &= \frac{1}{2} \left( -\frac{b_1+b_2}{b_2} -
    \frac{b_1+b_2}{z-1} \right) 
    \left( \frac{-b_2}{z-1} -1 \right), \\
  A_{10} &= \frac{1}{2} \cdot
    \left( \frac{b_2(z-(a_1+a_2))}{(b_1+b_2)(z-1)} -1
    \right)^2, \\
  A_{11} &= \frac{1}{2} \left( \frac{-b_1(z-(a_1+a_2))}{(b_1+b_2)(z-1)}
    + \frac{b_1}{b_2} \right) \left(
    \frac{b_2(z-(a_1+a_2))}{(b_1+b_2)(z-1)} - 1 \right), \\
  A_{12} &= \frac{1}{2} \left( -\frac{b_1}{b_2} -
    \frac{a_1+a_2-1}{a_1+a_2-(b_1+b_2)-1} \right) 
    \left( \frac{-b_2}{a_1+a_2-(b_1+b_2)-1} -1 \right), \\
  A_{13} &= \frac{1}{2} \cdot \left( c_2-1 \right)^2, \\
  A_{14} &= \frac{1}{2} \cdot \left( \frac{b_1}{b_2}-c_1 \right)
    \left( c_2-1 \right), \\
  A_{15} &= \frac{1}{2} \left( -\frac{a_1+a_2-1}{a_1+a_2-(b_1+b_2)-1}
    - \frac{b_1}{b_2} \right) \left(
    \frac{-b_2}{a_1+a_2-(b_1+b_2)-1} - 1 \right), \\
  A_{16} &= \frac{1}{2} \left( 1-(c_1+c_2)-\frac{a_1+a_2-1}{z-1}
    \right)  \left( c_2 - \frac{z-a_2}{z-1} \right), \\
\end{align*}

\begin{align*}
  A_{17} &= \frac{1-a_2}{z-1}
    \left( 1-(c_1+c_2) - \frac{a_1+a_2-1}{z-1} \right). \\
\end{align*}
We obtain
$$
\int\limits_{f \in R_1} +
\int\limits_{f \in R_2} =
  \begin{cases}
    0 & \textup{if } 1 < z < w, \\
    A_1 + A_3 & \textup{if } w \le z \le \frac{a_2-b_2}{a_2}, \\
    A_2 + A_3 & \textup{if } \frac{a_2-b_2}{a_2} < z < + \infty, \\
  \end{cases}
$$

$$
\int\limits_{f \in R_3} +
\int\limits_{f \in R_4} =
  \begin{cases}
    0 & \textup{if } 1 < z < a_1-c_1, \\
    A_4 + A_6 & \textup{if } a_1-c_1 \le z \le \frac{b_1-c_1}{b_1}, \\
    A_4 - A_5 + A_6 & \textup{if } \frac{b_1-c_1}{b_1} \le z < + \infty, \\
  \end{cases}
$$

$$
\int\limits_{f \in R_5} =
  \begin{cases}
    0 & \textup{if } 1 < z < a_1+a_2-(b_1+b_2) \\
    & \textup{\ \ or } a_1+a_2-b_1-1 > 0, \\
    A_7 - A_8 - A_9 & \textup{if } a_1+a_2-(b_1+b_2) 
    \le z \le 1-b_2 \\
    & \textup{\ \ and } a_1+a_2-b_1-1 \le 0, \\
    A_7 - A_8 & \textup{if } 1-b_2 \le z < + \infty \\
    & \textup{\ \ and } a_1+a_2-b_1-1 \le 0, \\
  \end{cases}
$$

$$
\int\limits_{f \in R_6} =
  \begin{cases}
    0 & \textup{if } 1 < z < a_1+a_2-(b_1+b_2), \\
    A_{10} - A_{11} - A_{12} & \textup{if }
    a_1+a_2-(b_1+b_2) \le z \le
    \frac{a_1+a_2-(c_1+c_2)}{1-(c_1+c_2)} \\
    & \textup{\ \ and } a_1+a_2-b_1-1 \le 0, \\
    A_{10} - A_{11} & \textup{if }
    a_1+a_2-(b_1+b_2) \le z \le
    \frac{a_1+a_2-(c_1+c_2)}{1-(c_1+c_2)} \\
    & \textup{\ \ and } a_1+a_2-b_1-1 > 0, \\
    A_{13} - A_{14} - A_{15} + A_{16} + A_{17} & \textup{if }
    \frac{a_1+a_2-(c_1+c_2)}{1-(c_1+c_2)} \le z <
    +\infty \\
    & \textup{\ \ and } a_1+a_2-b_1-1 \le 0, \\
    A_{13} - A_{14} + A_{16} + A_{17} & \textup{if }
     \frac{a_1+a_2-(c_1+c_2)}{1-(c_1+c_2)} \le z <
    +\infty \\
    & \textup{\ \ and } a_1+a_2-b_1-1 > 0. \\
  \end{cases}
$$

This leads to the lower bound for $P^{T_3}(z)$ which we
wanted, namely
\begin{equation*}
  P^{T_3}(z) \ge \frac{1}{\area(T_3)} \sum_{j=1}^{6}{
  \int\limits_{f \in R_j}{\boldsymbol{1}\{\bar{t} \le z\}df}}.
\end{equation*}
As for quadrilaterals, we did not succeed in showing
algebraically the convergence of this lower bound to $1$
when $w$ converges to $1$.
However, simulations with discretized parameters
$a_1 \in \{3, 4, \dots, 1\,000\,000\}$,
$a_2 \in \{0{.}001, 0.00{2}, \dots, 0{.}999\}$, and
$b_1 \in \{0{.}001, 0{.}002, \dots, 0{.}999\}$ (where $b_1 <
\frac{a_2}{a_1 + a_2 - 1}$)
such that $w$ converges to $1$ suggest that this is the
case.

\bibliographystyle{plain}
\bibliography{probability}

\end{document}